\documentclass[envcountsect]{svmult}  

\usepackage{makeidx}         
\usepackage{graphicx}        
\usepackage{multicol}        
\usepackage[bottom]{footmisc}

\usepackage{amsmath}
\usepackage{amsfonts}

\usepackage{epstopdf}	

\usepackage{dsfont}
\newcommand{\Ind}{\mathds{1}}

\newcommand{\R}{{\mathbb {R}}}
\newcommand{\N}{{\mathbb {N}}}
\newcommand{\Z}{\mathbb {Z}}
\newcommand{\C}{{\mathbb{C}}}
\newcommand{\Esp}{{\mathbb{E}}}
\newcommand{\sphere}{\mathbb{S}}

\newcommand{\tr}{\textmd{trace}\,}

\spnewtheorem*{conjectur}{Conjecture}{\bf}{\it} 
\spnewtheorem*{rk}{Remark}{\bf}{} 
\spnewtheorem*{rkss}{Remarks}{\bf}{} 
\newenvironment{rks}%
{\begin{rkss}\begin{enumerate}}%
{\end{enumerate}\end{rkss}}%

\begin{document}

\title*{The Hypergroup Property and Representation of Markov Kernels}
\author{Dominique Bakry\inst{1}\and Nolwen Huet\inst{2}}
\institute{%
IUF and Institut de Math\'ematiques,
Universit\'e de Toulouse and CNRS\\ 
118, route de Narbonne, 31400 Toulouse, France \\
\textit{dominique.bakry@math.univ-toulouse.fr} \and
Institut de Math\'ematiques, Universit\'e de Toulouse and CNRS\\ 
118, route de Narbonne, 31400 Toulouse, France
\textit{nolwen.huet@math.univ-toulouse.fr}
}
\maketitle

\begin{abstract}
For a given orthonormal basis $(f_n)$ on a probability measure space, we want to describe   all Markov operators which have the $f_n$ as eigenvectors. 
We introduce for that  what we call the  \emph{hypergroup property}. We study this property in three different cases.

On finite sets, this property appears as the dual of the GKS property linked with correlation inequalities in statistical mechanics. The representation theory of groups provides generic examples where these two properties are satisfied, although this group structure is not necessary in general.

The \emph{hypergroup property} also holds for Sturm--Liouville bases associated with log-concave symmetric measures on a compact interval, as stated in Achour--Trim\`eche's theorem. We give some criteria to relax   this symmetry condition in view of extensions to a more general context.

In the case of Jacobi polynomials with non-symmetric parameters, the  \emph{hypergroup property} is nothing else than Gasper's theorem. The  proof we present is based on a natural interpretation of these polynomials as harmonic functions and is related to analysis on spheres. The proof relies on the representation of the polynomials     as the moments of a complex variable.
\end{abstract}

\section*{Contents}
\setcounter{minitocdepth}{2}
\dominitoc

\section{Introduction}\label{intro}
In a number of situations, Markov operators appear to be a wonderful tool 
to provide useful information on a given measured space. 
Let us for example mention heat kernel methods to prove functional inequalities
like Sobolev or Log-Sobolev 
inequalities, or Cauchy kernels to prove boundness 
results on Riesz transforms in $L^p$. Heat kernels are widely used in 
Riemannian geometry and statistical mechanics, while Poisson, Cauchy 
and other kernels had been proved useful in other contexts related to 
classical harmonic analysis (see \cite{bakry-emery2-85,bakry-riesz2,bakry-riesz1,bakry-riesz3,beckner3,beckner-fourier,davies,Gross,gross2,LivreLogSob,scheffer} for example to see the 
action of different families of semigroups in various contexts).

It seems therefore interesting to describe all Markov kernels associated with a given structure. 
In what follows, we shall consider a probability space 
$(E, {\cal E}, \mu)$ on which is given an orthonormal $L^{2}(\mu)$ basis
$\mathcal{F}=(f_{0},\ldots,f_{n},\ldots)$, where we impose $f_{0}= 1$. 
Such a basis shall be called a unitary orthonormal basis (UOB in short).  
In general our basis $\mathcal{F}$ shall be real, but we do not exclude to 
consider complex bases.

We may then try to describe the  Markov operators defined from a family of probability measures $k(x,dy)$ by
$$f\mapsto K(f) = \int _{E}f(y)\ k(x,dy)$$ which are symmetric in $L^{2}(\mu)$ and have the functions $f_{n}$ as eigenvectors.  In other words,  we want to define the linear operator  $K$ from
$$K(f_{n})= \lambda_{n} f_{n},$$
 and try to describe for which sequences $(\lambda_{n})$ this operator is a Markov kernel. 
 We shall call these sequences Markov sequences  (MS's in short) associated with the UOB $\mathcal{F}$.
 
 For a general basis $\mathcal{F}$, this is quite impossible. But many bases which appear in natural
 examples have a special property, which we call the hypergroup property, under which one is able to 
 describe all Markov sequences  associated with the UOB $\mathcal{F}$.
 
 This expository paper is not intended to be a complete account of the general theory of hypergroups, 
 for which we may for example refer to the complete treatise 
 \cite{BloomHeyer95}. In fact, we just extracted 
 from this theory what is useful for our purpose.  More precisely, we 
 concentrated on the fundamental aspect which we are interested in, 
 that is the possibility of describing all Markov kernels associated 
 with our basis $\mathcal{F}$.

 The paper is organized as follows. 
 
 In the first part, we present the case of finite sets, where the hypergroup property appears as the 
 dual property of a more natural condition on the basis $\mathcal{F}$, namely the positivity of the multiplication 
 coefficients.  This property is called the GKS property in 
 \cite{bakry-echerbault-GKS} because of its links with 
 some famous correlation inequalities in statistical mechanics, and 
 we keep this notation.  These correlation inequalities 
 did in fact motivate our efforts in this direction (see Paragraph 
 \ref{GKS2}). Many 
 examples come from the representation theory of groups, but we 
 propose a systematic exposition.  The hypergroup property 
 provides a convolution operation on the set of probability measures on 
 the space, and the Markov kernels may be represented as the 
 convolution with some given measure.  In this situation, the hypergroup 
 property appears as a special property of an orthogonal matrix. We shall see that in fact there are 
 many situations where no group structure holds and where nevertheless 
 the hypergroup property  holds.
 
 The second part is devoted to the presentation of Achour--Trim\`eche's theorem, which states this property for the 
 basis of eigenvectors of a Sturm--Liouville operator with Numen boundary conditions, associated with a 
 log-concave symmetric measure on a compact interval. The original 
 Achour--Trim\`eche's theorem was not stated exactly in the same way (see \cite{AchourTrimeche}), but his 
 argument carries over very easily to our context.  We give a complete proof of this result, since to 
 our knowledge this proof was never published.   We tried to relax the 
 symmetry condition on the measure, and provided various rather 
 technical extensions of the theorem. But we did not succeed to extend 
 Achour--Trim\`eche's theorem to a wider class of measures which would include 
 the case of Gasper's theorem on Jacobi polynomials studied in the 
 next chapter. The real motivation of this effort is that in 
 general, log-concave measures on $\R$ or on an interval serve as a 
 baby model  in Riemannian geometry for manifolds with non-negative 
 Ricci curvature. Unfortunately, the symmetry condition does not seem 
 to have any natural interpretation.
 
 In the last part, we present Gasper's theorem, which states hypergroup property for the Jacobi polynomials.
 We present a proof which relies on geometric considerations on the spheres when the parameters are 
 integers, and which easily extends to the general case. We follow 
 Koornwinder's proof of the result, and give a natural interpretation 
 of Koornwinder's formula (Lemma \ref{lem2}) which represents those 
 polynomials as the moments of some complex random variable. We found after 
 the redaction of this part  that the interpretation of Jacobi polynomials as harmonic functions was already known 
 from specialists (see \cite{braak-meul,Koorn71,Koorn73}) but it seems that it was not directly 
 used to prove this integral representation formula. We hope that 
 this simple interpretation may provide other examples for similar 
 representation in other contexts.
 
 \section{The finite case}
 
 In this section, we restrict ourselves to the case of finite sets, 
 since in this context most of the ideas underlying the general 
 setting are present, and we so avoid the analytic complexity of 
 the more general cases that we shall study later on.

\subsection{The GKS property}

In what follows, we assume that our space is a finite set
$$E= 
\{x_{0}, \ldots,x_{n}\},$$ endowed 
with a probability measure 
$$\mu= (\mu(x_{0}), \ldots, \mu(x_{n})).$$  We denote by $L^{2}(\mu)$ the space of real functions on $E$, and we assume the existence of a real basis
$$\mathcal{F}= (f_{0}, \ldots, f_{n})$$ with $f_{0}=1$. We suppose here that 
for any $x\in E$, $\mu(x)>0$.

We shall write $\langle f\rangle$ for $\int f\ d\mu$ and 
$\langle f,g\rangle$ for $\int fg\ d\mu$.

The algebra structure of the set of functions is reflected in the 
multiplication tensor $(a_{ijk})$ for which
$$f_{i}f_{j} = \sum_{k} a_{ijk}f_{k}.$$

We therefore have
$$a_{ijk}= \langle f_{i}f_{j}f_{k}\rangle,$$
and we see that the tensor $(a_{ijk})$ is symmetric in $(ijk)$. It 
has also another property which reflects the fact that the 
multiplication is associative.

\begin{definition} We shall say that $\mathcal{F}$ has the property GKS if all 
coefficients $(a_{ijk})$ are non-negative.
\end{definition}

This notation comes from the GKS inequality in statistical mechanics 
that we shall describe at the end of this section.

Observe that $a_{ij0}= \delta_{ij}$.

Many natural bases $\mathcal{F}$ share this property. For example, consider 
the hypercube 
$E=\{-1,1\}^N$, with  the uniform measure on it.
 Let $\omega_{i}$ denote the $i$-th coordinates
 $$\omega= (\omega_{1}, \ldots,\omega_{n}) \mapsto \omega_{i},$$
 and, for $A\subset \{1, \ldots, N\}$
 $$\omega_{A}= \prod_{i\in A} \omega_{i}, ~\omega_{\emptyset}=1.$$
 Then, 
 $$\mathcal{F}= \{\omega\mapsto\omega_{A}, A\subset \{1, \ldots, N\}\}$$ is a UOB of 
 $(E,\mu)$. Since 
 $$\omega_{A}\omega_{B}= \omega_{A\Delta B},$$ it has the GKS property.
 
 (We shall see later that this is a special case of a generic 
 situation in finite groups).
 
 Although we are here mainly interested in the case of a real 
 basis, there are many natural complex GKS bases, issued in general
 from the representation theory of finite groups (see Paragraph 
 \ref{groupes} later). If the basis is complex, we shall still require 
 that the multiplication coefficients are non-negative real numbers, 
 which means that 
 $$a_{ijk}= \langle f_{i}f_{j}\bar f_{k}\rangle\geq 0,$$
for any $(i,j,k)$.

In what follows, we only consider real GKS bases, although the next 
result remains probably true in the complex setting.

\begin{proposition} \label{prop:GKS-x_0} If a UOB has the GKS property, then there exists a unique 
point $x_{0}$ on which every $f_{i}(x_{0})$ is maximal.  Moreover,  
for any $i$ and any $x$,
$\left|f_{i}\right|(x) \leq f_{i}(x_{0})$, and  at 
this point $x_{0}$, $\mu(x_{0})$ is minimal.

\end{proposition}

\begin{proof} Let us say that a function $f: E \mapsto \R$ is GKS if for any 
$i= 0, \ldots, n$, $\langle f f_{i}\rangle \geq 0$. In other words,  $f$ is written with non-negative coefficients  in the basis $\mathcal{F}$.

We shall say that a set $K\subset E$ is GKS if $\Ind_{K} $ is a GKS 
function. We shall say that a point $x$ is GKS if $\{x\}$ is a GKS set.

We shall see that there is only one  GKS point.

Remark first that the sum of two GKS functions is GKS and that,  thanks to the GKS property of $\mathcal{F}$, the 
product of two GKS functions is GKS. Moreover a limit of GKS functions 
is GKS. 
Observe also that a GKS function has always a non-negative integral 
with respect to~$\mu$ since $f_0=1$.
  
Let us consider a non-zero GKS function $f$ and consider 
$m= \max_{x\in  E}\left|f\right|$.   We see first that $\{f= m\}\neq 
\emptyset$.  

For this, assume the contrary, that is that $f=-m$ on 
$\left|f\right|=m$. Since $m>0$, we see that
 ${f^{2p+1}}/{m^{2p+1}}$ is a GKS function, and converges to 
 $-\Ind_{\{f=-m\}}$. Since a GKS function has a non-negative integral, 
 this is impossible.  Using the same argument, we see that
 $$\frac{1}{2}\left(\lim_{n}\frac{f^{2n}}{m^{2n}}+ 
 \lim_{n}\frac{f^{2n+1}}{m^{2n+1}}\right)= 
 \Ind_{\{f=m\}}$$ is a GKS function. 
 
Therefore, the set  
  $\{f=m\}$ is a GKS set, and there are non-trivial GKS sets.
 
Moreover, for any GKS function, \label{max_f_gks_set}
 the set $\{f= \max 
 \left|f\right|\}$ is GKS. 
 
 Let $E_{1}$ be a nonempty GKS set, minimal for the inclusion. Then, 
 for any GKS function $f$, $g=\Ind_{E_{1}}f$ is GKS. If $g$ is not $0$, 
 then its maximum is attained on a subset $E_{2}$ of $E_{1}$ which 
 is again GKS. Since $E_{1}$ is minimal, we have $E_{2}=E_{1}$.
 
 Therefore, for any GKS function, its restriction $g$ to $E_{1}$ is either  $0$ on $E_{1}$ or 
 constant (and equal to the  maximum of $g$). In any case, 
 $f$ is constant on $E_{1}$. Since this applies to every function 
 $f_{i}$ in $\mathcal{F}$, and since $\mathcal{F}$ is a basis, every function is constant 
 on $E_{1}$ and therefore $E_{1}$ is reduced to a single point 
 $\{x_{0}\}$.
 
 The same proof shows that any GKS set contains a GKS point.
 
 For a GKS point $x_{0}$, 
 $f_{i}(x_{0}) ={\langle f_{i} \Ind_{x_{0}}\rangle}/{\mu(x_{0})}\geq 
 0$.

Then, consider two distinct points $x_{0}$ and $x_{1}$, 
  and write
 $$\frac{\Ind_{x_{0}}}{\mu(x_{0}) }= \sum_{k} f_{k}(x_{0}) 
 f_{k},~\frac{\Ind_{x_{1}}}{\mu(x_{1}) }= \sum_{k} f_{k}(x_{1}) f_{k}.$$ 
 Writing the product, we see that
 $$0= 
\frac{\Ind_{x_{0}}\Ind_{x_{1}}}{\mu(x_{0})\mu(x_{1})}= \sum_{ijk} 
 a_{ijk}f_{i}(x_{0})f_{j}(x_{1})f_{k},$$
 with the multiplication coefficients $a_{ijk}$.
 
 So we see that for any pair of distinct points, and for any 
 $k$,
 $$\sum_{ij} f_{i}(x_{0})f_{j}(x_{1})a_{ijk} =0.$$
 
 Suppose then that $x_{0}$ and 
 $x_{1}$ are GKS points. In the previous sum, all coefficients are non-negative.  
 Therefore, for any $(i,j,k)$
 $$a_{ijk} f_{i}(x_{0})f_{j}(x_{1})=0.$$
  If we apply that with $i=0$ and $j=k$, we see that 
  $f_{j}(x_{1})=0$, for any $i$. This is impossible since then any 
  function would take the value $0$ in $x_{1}$. So there is a unique GKS point.
  
  Let $x_{0}$ be this unique GKS point.  Any GKS set contains 
  $x_{0}$. Since for any GKS function, the set where $f$ is maximum 
  is GKS, any GKS function attains its maximum at $x_{0}$.

  It remains to show that $
\mu$ is minimal at $x_{0}$.   For this, 
  observe that for any point $x$, the function
  $$f_{x}= \mu(x) \Ind_{x_{0}}- \mu(x_{0}) \Ind_{x}$$ is GKS (this comes 
  from the fact that each $f_{i}$ is maximal at $x_{0}$).
   Therefore, the maximum value of $\left|f_{x}\right|$ is attained 
   in $x_{0}$, which gives the result.\qed
\end{proof}

\subsection{Orthogonal matrix representation}

Consider the matrix  $$(O_{ij})= 
\left(\sqrt{\mu(x_{i})}f_{j}(x_{i})\right),$$
 we see easily that the matrix $(O_{ij})$ is a $(n+1)\times (n+1)$ orthogonal matrix 
 with positive first column.
 Conversely, any such matrix may be associated with a UOB 
 on a finite set with measure $\mu$ given by
 $$\mu(x_{i})= O_{i0}^{2}.$$
  Therefore, there is a one to one correspondence between the set of 
  orthogonal matrices with positive first column, and the set of 
  finite probability spaces, whose probability has everywhere positive weight, endowed with a UOB. (In fact, this is not 
  completely true, since we would not distinguish between bases given 
  in different orders, provided that the first element is $1$, which 
  identifies the set of UOBs with a quotient of a the set of 
  orthogonal matrices through a permutation of rows and columns.)
  
  The GKS property may be translated into the following property on such 
  an orthogonal matrix:
  
  \begin{equation}\label{GKSmatrix}\forall j,k,l, ~ \sum_{i} 
  \frac{O_{ij}O_{ik}O_{il}}{O_{i0}}\geq 0.\end{equation}
  
  The transposed of an orthogonal matrix is orthogonal, and we just 
  saw that an orthogonal matrix which has the GKS property also has a
  non-negative row (corresponding to the row where the first column is minimal according to Proposition \ref{prop:GKS-x_0}). We may of course rearrange the labelling of the 
  points in such a way that this row is the first one. Then, the 
  situation is completely symmetric.

 We shall then consider the squares of terms in the first 
row as a probability measure on the dual set $\{0,1,\ldots,n\}$: 
$$\nu(i)= \mu(x_{0}) f_{i}^{2}(x_{0}).$$

Thanks to the fact that the functions $f_{i}$ are 
maximal at $x_{0}$ and that this maximum must be larger than $1$ 
(since $\int f_{i}^{2}\ d\mu=1$), the dual measure is also minimum at 
$0$.

As an application, we have the following.

\begin{proposition}\label{GKSUniforme}
If a real GKS basis exists for the uniform measure on some finite set $E$, then the cardinal  of $E$ must be $2^{k}$ for some $k$, and this basis is the canonical basis $(\omega_{A})$ of the characters of the group $(\Z/2\Z)^{k}$.

\end{proposition}

To see this, we first observe that the dual measure is uniform too.   
In fact, for the matrix $(O_{ij})$, we have  $O_{00}= 
{1}/{\sqrt{n+1}}$, where $n+1$ is the number of points in the space, 
and since the first row is positive and has minimum value 
${1}/{\sqrt{n+1}}$, it must be constant since the sum of the 
squares of its coefficients is $1$. 

Now, if we multiply the matrix $O$ by $\sqrt{n+1}$, then we see that in 
each column, the maximum value of the coefficients is attained on the 
first row and is equal to 1.
 Since the sum of 
all the squares of the coefficients in a given column must add to 
$n+1$, this shows that in any column, the coefficients must take only the
values $\pm 1$.

Those matrices (with entries $\pm 1$ and orthogonal lines) are called Hadamard matrices (cf \cite{hadamard,paley-hadamard,Kharaghani-Tayfeh-Hadamard}). 
If  $n+1= 2^{k}$ for some $k$,  such matrices are given by the basis  
 $(\omega_{A})$ on $\{-1,1\}^k$ and satisfy the GKS property. 
 It is known that the order of a Hadamard matrix must be $1$, $2$, or $4k$, 
 and it is an open problem to find such matrices for all $k$ 
 (the lowest $k$ for which no Hadamard matrix of order $4k$ is known is 
 $k=167$ since \cite{Kharaghani-Tayfeh-Hadamard}).
 Nevertheless, the following  proposition will prove the result of Proposition  \ref{GKSUniforme}.
 \begin{proposition}
 If a Hadamard matrix has the GKS property, then it must be of order  
 $2^{k}$ for some $k$,  and, up to permutation,  it is  the matrix of the 
 canonical basis $(\omega_{A})$ 
   of the group $\{-1,1\}^k$.
 
 \end{proposition} 
 
 \begin{proof}  The case of a set of size 2 is trivial, and we therefore assume that the size of the matrix is at least 3. To fix the idea, consider  a matrix $M$ of order $n+1$ with entries 
 $\pm 1$, with orthogonal columns.  Call $f_{i}$ the column vectors 
 and $\hat f_{i}= 
{f_{i}}/{\sqrt{n+1}}$ the normalized ones, so that $M=\big(f_i(x_j)\big)_{0\leq i,j\leq n}$.  As usual, let us denote by $\langle 
f\rangle$ the mean value of a function $f$ with respect to the 
normalized uniform measure. We suppose, which is possible up to reordering,  that all the entries of the first line and the first 
column are $+1$. The GKS property says that 
$$\langle f_{i}f_{j}f_{k}\rangle \geq 0,$$ for any $(i,j,k)$.

First, since $\langle f_{i}\rangle =0$ for any $i\geq 1$, there must be as 
many $1$ and $-1$ in each column, and therefore  $n+1$ is even. Let $A$ be the set of points where 
$f_1=1$.
Let $i\geq 2$, and let $q$ be the number of points in $x\in A$ such that 
$f_{i}(x)=1$, and $r$ be the number of points in $A^{c}$ where 
$f_{i}(x)=1$. Writing 
$\langle f_{i}\rangle =0$ and $\langle f_{1}f_{i}\rangle =0$, we get 
$$q+r=p, ~q=r,$$ which  shows that $p$ is even and also that,  when $i\geq 2$,
$$\langle  \Ind_{A} f_{i}\rangle =0.$$ This is the generic 
argument which shows that Hadamard matrices have order $4k$.
We shall now make use of the GKS property.
Write $$\Ind_{A}f_{i}= \sum_{l} a_{il}f_{l}.$$
We have $$a_{il}= \langle \Ind_{A} \hat f_{i} \hat f_{l}\rangle= 
\frac{1}{n+1}\langle \Ind_{A} f_{i} f_{l}\rangle ,$$ and therefore the matrix 
$A=(a_{il})$ is symmetric. As it is the matrix of a projector, 
its eigenvalue are $0$ or $1$. The dimension of the eigenspace 
associated to $1$ is $p$, since the eigenspace is generated by 
$(\Ind_{\{x\}}, x\in A)$.
Also, $A$ is a GKS set because it is the set where $f_{1}$ 
attains it's maximum (cf proof of Proposition \ref{prop:GKS-x_0}, Page \pageref{max_f_gks_set}). This implies  that all the entries of $A$ are not 
negative.
Moreover, if we look at the values of the functions at $x_0$, it holds $\sum_{j} a_{ij}=1$. 

Therefore, the matrix $A$ is Markovian. There are no transitory points 
since $A$ is symmetric.  The number of recurrence classes for such a matrix is the 
multiplicity of $1$ as eigenvector, here $p$, so there are exactly $p$ 
recurrence classes, and no recurrence class is reduced to a single 
point, since $a_{ii}=1/2<1$.
So every recurrence class has exactly two points.  For example, $\{0,1\}$  form a recurrence class. On each line of the matrix, there are exactly two places where $a_{ij}\neq 0$. The values of those entries are then $1/2$, since $a_{ii}= 1/2$. If we choose 
two distinct indices $i$ and $j$ in two different recurrent classes, 
then  $a_{ij}=0$. This means that $\langle \Ind_{A} f_{i}f_{j}\rangle 
=0$. 

Choose now one index in every recurrence class (say the even indices to
fix the ideas, which is possible up to reordering of the columns).
Then, the functions $g_{i}=\Ind_{A}f_{2i} $ form an orthogonal Hadamard GKS
matrix of order $p=(n+1)/2$, and we may now use induction to see that the order must be $2^{k}$ for some $k$.

To see that the unique basis such basis in dimension $2^{k}$ is given 
by the canonical basis, it is enough to observe that if $i$ and $\sigma(i)$ 
are in the same recurrence class for the matrix $A$, then we have 
\[
\Ind_{A}f_{i}=\frac{1}{2} (f_{i}+ f_{\sigma(i)}) =\Ind_{A} f_{\sigma(i)},
\]
and also 
\[
\Ind_{A^{c}}f_{i}= -\Ind_{A^{c}}f_{\sigma(i)}.  
\]
Then an easy induction leads to the result.\qed
\end{proof}

Nevertheless, unlike real bases, we shall see in Paragraph \ref{groupes} that there 
always 
exists a complex GKS basis on any finite set with uniform measure.

\label{exdim2}
On two points, an easy 
computation shows that a two dimensional orthogonal matrix 
having the GKS property must be
     $$\begin{pmatrix} \cos \theta&\sin \theta\\ \sin 
     \theta&-\cos\theta\end{pmatrix}$$
      with $\theta\in [\pi/4, \pi/2)$. In fact, in dimension $2$,  
      given the measure (giving two distinct positive masses on the 
      two points),
      there are exactly two UOB, and only one such that the unique non-constant function is maximal on the point with minimal mass (a necessary condition
      to have the GKS property as we saw). 
     In this situation, any GKS matrix is symmetric, and the set of 
orthogonal matrices having the GKS property is connected.
     
  On three points, the situation is more complicated. We saw for example 
  that  
      there are no real GKS basis when the measure is uniform. The set of real
      GKS UOBs on three points  is connected, and one may see that the 
      maximum value of $\mu(x_{0})$ for which there exists a real GKS 
      basis is $\mu(x_{0})=1/4$, and in this case the probability 
      measure is $(1/4, 1/4, 1/2)$ and the unique GKS basis is obtained 
      taking the real parts of the characters in the group $\Z/4\Z$ (we 
      shall see later in Paragraph \ref{groupes} how to associate 
      complex or real GKS bases with any finite group).
      There is also a complex GKS basis with the uniform measure (the 
      characters of the group $\Z/3\Z$).

\subsection{The hypergroup property}

As we saw before, when we have the GKS property, the situation is completely symmetric and we may consider the 
dual property.   This is the hypergroup property.

\begin{definition} We shall say that the UOB $\mathcal{F}$ of $L^{2}(\mu)$ satisfies the 
hypergroup property (HGP in short) at point $x_{0}$ if
 for any $(x,y,z) \in E^{3}$,
$$K(x,y,z)= \sum_{i}\frac{ f_{i}(x)f_i(y) f_{i}(z)}{f_{i}(x_{0})} \geq 0.$$

Of course, this supposes that for any $i$, $f_{i}(x_{0})\neq 0$.

 \end{definition}
 
 We do not require the GKS property to hold in the definition of the HGP property. We shall see later (at the end of Paragraph \ref{groupes}, Page \pageref{exdim3}) that we can find bases with HGP property and not GKS property, and the reverse. On the other hand, we shall see that in groups, the natural basis always share both property.
 
 \begin{rk}
 Since we may always change for any $i\geq1$  $f_{i}\in \mathcal{F}$ into $-f_{i}$, and 
 still get a UOB,  and that this operation does not change the 
 hypergroup property, we see that we may always assume that 
 $f_{i}(x_{0})>0$.
 \end{rk}

 To see the duality with the previous situation, let us enumerate the 
 points in 
 $E$ starting from $x_{0}$, and recall our orthogonal matrix $O$ with
 $O_{ij}= \sqrt{\mu(x_{i})}f_{j}(x_{i})$. 
 
Then, under the GKS condition, $O$ has  non-negative first line and first 
column. Recall that the GKS property may be written as 
$$\forall j,k,l, ~\sum_{i} \frac{O_{ij}O_{ik}O_{il}}{O_{i0}}\geq 0.$$
while the HGP property writes
 $$\forall j,k,l, ~ \sum_{i} \frac{O_{ji}O_{ki}O_{li}}{O_{0i}}\geq 0,$$

 Notice also that if both properties occur, then the point $x_{0}$ 
 must be the unique point where all functions $f_{i}$ are non-negative, and where all $f_i$ and $\left|f_i\right|$ are maximal.
 
 From the symmetry of the situation, we may consider the functions 
 $$g_{x}(i)= \frac{O_{xi}}{O_{0i}}$$ to be a UOB on the set 
 $\{0,\ldots,n\}$ endowed with the measure $\nu(i)= O_{0i}^{2}$.
 
 From this we deduce that  a basis has the HGP property if and only 
 if  this new 
 basis has the GKS property, where of course the point $0$ plays the 
 role of the point $x_{0}$ in the previous paragraph. Therefore, the 
 $g_{x}(i)$ are maximal at $i=0$, and moreover $0_{i0}$ is minimal at 
 $i=0$. This means that, for the HGP property also, one has
 
 \begin{proposition}
 If the UOB $\mathcal{F}$ has the HGP property at point $x_0$, then
 $$\forall i,x, ~\left|f_{i}(x)\right|\leq f_{i}(x_{0}); ~
  ~ \forall x, ~\mu(x)\geq \mu(x_{0}).$$
 (Recall that we assume that $f_i(x_0)>0$.)  
 \end{proposition}
 
  We may reformulate the HGP property in the following way, which 
  we shall use later in a different context, since it is more 
  tractable.
  
  \begin{proposition}\label{RepresentationIntegrale}
  The UOB $\mathcal{F}$ has the HGP property if and only if there is a 
  probability kernel $k(x,y,dz)$  such that, for any $i= 0,\ldots,n$
  $$\frac{f_{i}(x)f_{i}(y) }{f_{i}(x_{0})}= \int \frac{f_{i}(z)}{f_{i}(x_{0})}\ k(x,y,dz).$$
 \end{proposition}
 
 \begin{proof}
 The proof is straightforward. If the hypergroup property holds, then 
 the kernel 
 $$k(x,y,dz)= \left(\sum_{i} 
 \frac{f_{i}(x)f_{i}(y)f_{i}(z)}{f_{i}(x_{0})}\right) 
 \mu(dz)$$ is a probability kernel satisfying our conditions.
 
 On the other hand, if such a probability kernel $k(x,y,dz)$ exists, 
 then writing $$k(x,y,dz) = K_{1}(x,y,z) \mu(dz)$$ and
 $$K_{1}(x,y,z) = \sum_{i} a_{i}(x,y) f_{i}(z),$$
 one sees that 
 $$a_{i}(x,y)= \frac{f_{i}(x)f_{i}(y) }{f_{i}(x_{0})}.$$\qed
\end{proof}
  The link with the Markov operators is the following. 
  \begin{definition} 
   A Markov operator is just an operator $K$ which satisfies $K(1)=1$ 
   and which preserves non-negative functions.

   Given a sequence $\lambda= (\lambda_{i}, ~{i=0,\ldots,n})$, we 
   define the associated linear operator $ K_{\lambda}$ by 
   $$K_{\lambda}(f_{i})= \lambda_{i}f_{i}.$$
   A Markov sequence (MS 
   in short) is a sequence $\lambda=(\lambda_{i})$ 
   such that the associated operator $K_{\lambda}$ is a Markov 
   operator.
   
   \end{definition}
   
 Remark that for any MS $\lambda$, one has $\lambda_{0}=1$.
   Remark also that 
     $$K_{\lambda}(f)(x) = \int f(y)\ k_{\lambda}(x,dy) ,$$
      where 
      $$k_{\lambda}(x,dy) = \left(\sum_{i} \lambda_{i} f_{i}(x)f_{i}(y)\right) \mu(dy).$$
      Therefore, the set of MS's is just the set of sequences $\lambda$ such that the matrices 
      $k_{\lambda}(x,y)$ are Markov matrices.
   
   The HGP property asserts that, for any $x$, the sequence
   $$\lambda(x)= \left(\frac{f_{i}(x)}{f_{i}(x_{0})}\right)_{i=0,\ldots,n}$$ is a MS.
   
   It is quite standard to see that any eigenvalue $\lambda_{i}$ of a Markov 
   operator must satisfy $\left|\lambda_{i}\right|\leq 1$.
      The set of Markov sequences is a convex compact set, which is 
   stable under pointwise multiplication.
   Now   the main interest of this property relies in the following theorem.
   
   \begin{theorem}\label{thmfond}
   If the basis $\mathcal{F}$ has the HGP property, then the sequences $$\left(\lambda_{i}(x)= \frac{f_{i}(x)}{f_{i}(x_{0})}\right)_i$$ are 
   the extremal points in the convex set of all Markov sequences.
   
     More precisely,  every Markov sequence may be written uniquely as
   $$\left(\lambda_{i}= \int \frac{f_{i}(x)}{f_{i}(x_{0})}\ d\nu(x)\right)_i$$ for 
   some probability measure $\nu$ on $E$. Conversely, every probability measure can be associated in the same way with a Markov sequence.
     \end{theorem}

   The main interest of this result is that there exist numerous 
   natural $L^{2}$ bases with the HGP property, as we shall see later.
   
   \begin{proof}  The representation formula is straightforward. Indeed, 
   writing $\nu(dy)= k_\lambda(x_{0},dy)$, one has 
   $$\lambda_{i} f_{i}(x_{0})= \int f_{i}(y)\ \nu(dy),$$ which gives 
   the representation.
   
  From this, it is easy to see that if the sequences $(\lambda_i(x))_i$ are Markov sequences, then they are 
  extremal. Indeed any representation 
  $$\frac{f_{i}(x)}{f_{i}(x_{0})}= \theta \lambda_{i}^{1}+ 
  (1-\theta)\lambda_{i}^{2}$$ with MS's $\lambda^{1}$ and $\lambda^{2}$ 
  leads to
  $$f_{i}(x) = \int f_{i}(y)\ \nu(dy),$$ with 
  $$\nu(dy)= \theta K_{\lambda^{1}}(x_{0},dy)+ (1-\theta) 
  K_{\lambda^{2}}(x_{0},dy).$$
  
  From this we deduce that for any function $f$
  $$f(x)= \int f(y)\ \nu(dy),$$ and therefore $\nu= \delta_{x}$, which 
  gives the extremality.\qed
\end{proof}

  \begin{rk}
       Remark that the representation formula is still true for any MS when the basis does not verify the HGP property. As we can always embed the convex set of Markov sequences in the $n$-dimensional affine space $H_{1}= \{(\lambda_{i}), ~\lambda_{0}=1\}$, the latter fact means that this set is actually contained in the $n$-simplex generated by the $n+1$ points $\lambda^{x/x_{0}}= 
  \big({f_{i}(x)}/{f_{i}(x_{0})}\big)_i$.
  
  Then, when the hypergroup property holds, the set 
  of Markov sequences  is a $n$-simplex and the 
  representation of a point in this set as affine combination of 
  extremal points is unique. 
 
  We may ask for which kind of $L^{2}$ basis on a finite space this 
  still happens. It is quite clear that the cardinal of the set of 
  extremal Markov sequences is finite. Indeed, 
  the set of Markov sequences is delimited by a finite number of 
  $(n-1)$-hyperplanes in $H_{1}$. Namely, for any pair $(x,y)$ of 
  points in $E$, one considers the half space defined by $\{\langle 
  \lambda, F^{x,y}\rangle \geq 0\}$, where $F^{x,y}= (f_{i}(x)f_{i}(y))$. 
  Then the set of Markov sequences is the intersection of all these 
  half spaces. Therefore, every extremal point  lies in the finite 
  set $E_{1}$ of
  possible intersections of $n$ hyperplanes $H^{x,y}=H_{1}\cap \{\langle 
  \lambda,F^{x,y}\rangle=0\}$. Now, consider any point $x_{0}$ such that  
  for any index $i$,
  $f_{i}(x_{0})\neq 0$. The point $\lambda^{x/x_{0}}= 
\big({f_{i}(x)}/{f_{i}(x_{0})}\big)_i$ belongs to $H^{x_{0},y}$ for any 
  $y\neq x$, thanks to the orthogonality relations of the basis. 
  Therefore, those points $\lambda^{x/x_{0}}$ belong to the set $E_{1}$. 
  When $x_{0}$ is fixed and $x$ varies in $E$, those points describe 
  a simplex $S_{x_{0}}$ for which we know that every Markov sequence belongs to 
  it. The hypergroup property holds at some point $x_{0}$ exactly 
  when no other 
   point in $E_{1}$ lie in the interior of $S_{x_{0}}$.
   
   On three points, one may check directly that the hypergroup 
   property holds at some point $x_{0}$ exactly when the set of 
   Markov sequences is a simplex (that means no other simplex is 
   possible than the simplices $S_{x},~ x\in E$). We may wonder if 
   this situation is general, that is if the hypergroup property is 
   equivalent to the fact that the set of Markov sequences is a 
   simplex.
 \end{rk}

   \subsection{Markov operators as convolutions}
   When the hypergroup property holds, we may introduce
   a convolution on the space of measures.

   Indeed, consider the kernel
      $$k(x,y,z)= \sum_{i}\frac{ f_{i}(x)f_i(y) 
      f_{i}(z)}{f_{i}(x_{0})}.$$ We observe that, for any $(x,y)$
      $$\int k(x,y,z)\ \mu(dz) =1,$$ and therefore the measures 
      $$\mu_{x,y}(dz) = k(x,y,z)\ \mu(dz)$$ are probability measures.
      
      We   may decide  that the convolution is defined from this 
      kernel by
	 $$\delta_{x}*\delta_{y}= \mu_{x,y},$$ and extending it to any 
	 measure by bilinearity.
    
   Moreover, we extend the convolution to functions by identifying a 
   function 
   $f$ with the measure $fd\mu$. This gives
   $$f*g(z)= \int f(x)g(y) k(x,y,z)\ d\mu(x)d\mu(y).$$
   
   Observe that \begin{equation}\label {eq1}f_{i}*f_{j}= \delta_{ij} 
   \frac{f_{i}}{f_{i}(x_{0})},\end{equation}
   
   and that this property again completely determines the convolution. 
   
   It is easy to verify that this convolution is commutative and that 
   $\delta_{x_{0}}*\nu= \nu$ for any $\nu$.
   Moreover,  if an operator $K$ satisfies $K(f_{i})= \lambda_{i} 
   f_{i}$, then
   $$K(f*g)= K(f)*g= f*K(g),$$
    as may be verified directly when $f= f_{i}$ and $g= f_{j}$ using 
    \eqref{eq1}.

    On the other hand, if $\nu$ is a probability measure, then the 
    operator $K_{\nu}(f)= f*\nu$ is a Markov operator which satisfies
    $$K_{\nu}(f_{i}) = \lambda_{i} f_{i},$$
     with
     $$\lambda_{i} = \frac{\int f_{i}\ d\nu}{f_{i}(x_{0})}.$$ This is straightforward using 
     \eqref{eq1} if we write $\nu(dx)= h(x)\mu(dx)$ and the 
     decomposition of $h$ along the basis $\mathcal{F}$.
    
    Therefore, if $K$ is a Markov operator, then we have
    $$K(f)= K(f*\delta_{x_{0}})= f*K(\delta_{x_{0}}).$$
     This representation is exactly the representation of Markov 
     sequences, with 
     $\nu= K(\delta_{x_{0}})$, and every Markov operator 
     $K_{\lambda}$ may be defined from
     $K_{\lambda}(f)= f*\nu$, for some probability measure $\nu$.

     \subsection{The case of finite groups}\label{groupes}
     
     Many natural examples of finite sets endowed with a probability 
     measure and a UOB which satisfies both GKS and HGP properties 
     come from finite groups.
     
     Since perhaps not every reader of these notes is familiar with 
    this setting,  let us summarize briefly the basic elements of the analysis on 
     groups. We refer to \cite{diaconis88} or \cite{Isaacs94} for 
     more details.
     
     Given a finite group $G$, one may consider linear 
     representations $\rho ~:~G\mapsto U(V)$, that is group homomorphisms between $G$ and 
     some $U(V)$, for some finite dimensional Hermitian space $V$ 
     (where $U(V)$ denotes the unitary group of $V$). Such a 
     representation is irreducible if there is no non-trivial proper subspace of 
     $V$ which is invariant under $\rho(G)$. Any representation may 
     be split into a sum of irreducible representations, acting on 
     orthogonal subspaces of $V$. Two representations 
     $(\rho_{1},V_{1})$ and $(\rho_{2},V_{2})$ are equivalent 
     if there exists a linear unitary isomorphism $h: V_{1}\mapsto V_{2}$ 
     such that $\rho_{1}(g)= h^{-1}\rho_{2}(g)h$ for any $g\in G$.
     There are only a finite number of non-equivalent irreducible 
     representations, that we denote $(\rho_{i},V_i)$, $i\in I=
     \{0,\ldots, n\}$.

     Let $\hat G$ the set of the equivalence classes of $G$ under the 
     conjugacy relation ($g_{1}$ is conjugate to  $g_{2}$ means $g_{2}= 
     g^{-1}g_{1}g$ for some $g\in G$). We endow $\hat G$ with the 
     probability $\nu$ which is the
     image measure of the uniform measure on $G$, which means that the measure of any class is proportional 
     to the number of points in this class. A function on $G$ which 
     is constant on conjugacy classes (that we call a class function) 
     can be seen as a function on $\hat G$. It is just a function which 
     is stable under conjugacy.

     For any irreducible representation $(\rho_{i}, V_{i})$, let us define the function $\chi_{i}$ on $G$ by $\chi_{i}(g)= \tr(\rho_{i}(g))$.
     This is a 
     class function, that is to say constant on any conjugacy class.  
     The function $\chi_{i}$ is called the character of the representation. By 
     convention, we take $\chi_{0}=1$, that is the trace of the  
     constant 
     representation into the space $V= \{\C\}$.
      
     \begin{proposition} The set $\{\chi_{i}, i\in I\}$ is a (complex) UOB for 
     $(\hat G, \nu)$. Moreover, it has the GKS and HGP properties.
     \end{proposition}
     
     \begin{proof} We shall not enter in the details here. We refer to any 
     introduction book on the representation theory of finite groups 
     for the first fact. We shall detail a bit more the HGP and GKS 
     properties, which are perhaps less standard.
     
     For the GKS property, for any pair of irreducible representations 
     $(\rho_{i},V_{i})$ and $(\rho_{j},V_{j})$, one may consider the 
     representation $\rho_{i}\otimes \rho_{j}$ in the tensor product 
     $V_{i}\otimes V_{j}$. If we split this  representation into 
     irreducible representations and take the trace, and if we notice 
     that $\tr (\rho_{1}\otimes \rho_{2})= 
     \tr(\rho_{1})\tr(\rho_{2})$, then we get that
     $$\chi_{i}\chi_{j}= \sum_{k}m_{ijk} \chi_{k},$$
     where $m_{ijk}$ is the number of times that the representation 
     $\rho_{k}$ appears in this decomposition. Here we may see that 
     not only the multiplication coefficients are non-negative, but they 
     are integers.
     
     We shall see next that this basis has the HGP property at the point 
     $x_{0}=e$ (which forms a conjugacy class by itself). For that, we require a bit more material.
     
     First define the convolution on the group 
     $G$ itself by
     $$\phi*\psi(g)= \frac{1}{\left|G\right|}\sum_{g'\in G} 
     \phi(gg'^{{-1}})\psi(g').$$
     
     The Fourier transform is defined on the set $I$ of irreducible 
     representation as
     $$\hat \phi(i)= \sum_{g\in G}\phi(g)\rho_{i}(g).$$
     (It takes values in the set of linear operators on $V_{i}$.)
     
     One has an inversion formula
     $$\phi(g)= \frac{1}{\left|G\right|}\sum_{i} 
     d_{i}\tr\left(\rho_{i}(g^{-1})\hat \phi(i)\right),$$ where 
     $d_{i}$ is the dimension of $V_{i}$ (the degree of the 
     representation).
     
     One has $$(\phi*\psi)\hat{}= \hat \phi \hat \psi,$$
     and 
     $$\hat\chi_{j}(i) = \delta_{ij} \frac{\left|G\right|}{d_{i}}.$$
     
     Now, the convolution of two class functions is again a class 
     function, as seen directly from the definition. 
     
     We want to show that this convolution is exactly the convolution 
     that we defined in the previous section from the HGP property, 
     that is 
     $$\chi_{i}*\chi_{j} = \delta_{ij}\frac{\chi_{i}.}{\chi_{i}(e)}.$$
     For that, we look at the Fourier transform and the result is 
     straightforward, since $\chi_{i}(e)=d_{i}$.
     
     This convolution is then the convolution defined from the 
     $\chi_{i}$, and we have
     $$\delta_{x}*\delta_{y}= k(x,y,z) d\mu(z),$$
     where 
     $$k(x,y,z) = \sum_{i} 
     \frac{\chi_{i}(x)\chi_{i}(y)\chi_{i}(z)}{\chi_{i}(e)}.$$
     
     Since by construction in this case the convolution of two 
     probability measures is a probability measure, the kernel 
     $k(x,y,z)$ is non-negative, which proves the HGP property.
     
     Observe that here the kernel $k(x,y,z)$ has a simple 
     interpretation. Given 3 classes $(x,y,z)$, then 
     $$k(x,y,z)=\frac{\left|G\right|}{\left|x\right|\left|y\right|}m(x,y,z),$$
      where $m(x,y,z)$ is, for any point $g\in z$, the number of ways 
      of writing $g=g_{1}g_{2}$ with $g_{1}\in x$ and $g_{2}\in y$, 
      this number being independent of the choice of $g\in z$.
     \qed
\end{proof}
     
     If we want to stick to real bases as we did before (and as we 
     shall do in the next chapters), we may restrict ourselves to 
     real groups (that is groups where $g$ and $g^{-1}$ are always in 
     the same class), or we may agglomerate the class of $g$ with the 
     class of $g^{-1}$. We get a new probability space, where the 
     functions $\Re(\chi_{i})$ form a UOB which again satisfies the 
     GKS and HGP properties.
     
     It is certainly worth noticing  that, unlike the 
     convolution on $G$ itself, the convolution on $\hat G$ is always
     commutative.
     
     Observe that taking the group ${\Z}/{n\Z}$, one gets a 
     complex GKS and HGP UOB on the set of finite points with the uniform measure 
     (with $f_{l}(x) = \exp(2i\pi lx)$), and that the unique real case 
     where the measure is uniform and is GKS (the hypercube) is 
     nothing else that the group $(\Z/2\Z)^{n}$.
     
     \begin{rk}
     Unlike what happens for finite groups, it is not true in general that a basis $\mathcal{F}$ which has the GKS  property has the dual property HGP. 
		 This is the case on two points spaces, since any orthogonal GKS 
	 matrix is symmetric (cf Page \pageref{exdim2}).
	  \label{exdim3}If we look at the sets with three points, one may construct examples of an orthogonal 
	  matrix having the  GKS property without the  HGP property (and 
	  conversely, of course). In fact, consider an orthogonal 
	  matrix $(O_{ij}), ~0\leq i,j\leq 2$, with positive first row 
	  and columns. $O_{00},O_{01}$ and $0_{10}$ determine 
	  entirely the first rows and columns, and then it is easy to see 
	  that there are only 2 orthogonal matrices with given 
	  $O_{00},O_{01},O_{10}$. Then, it is not hard (using a computer algebra 
	 program) to produce orthogonal matrices which have the GKS 
	 and not the HGP property, or which have neither, or both.
	\end{rk}

     \subsection{On the GKS inequalities}\label{GKS2}
     
     We conclude this section with some remarks on the correlation 
     inequalities in statistical mechanics.
     
     In this context, one is interested in the space of 
     configurations of some system. We have a set of positions $i\in 
     K$, $K$ being a finite set, and at each point $i\in K$ there is 
     some random variable $x_{i}$ with values in $E$, where $E$ is 
     some finite set, endowed with a probability measure $\mu$. One 
     is then interested in the set $E^{K}$ of configurations, which is 
     equipped with a measure $\mu_{H}$, where
     $$\mu_{H}(dx) = \exp (H(x)) \frac{\mu_{0}(dx)}{Z_{H}},$$
      where $\mu_{0}$ is the product measure $\mu^{\otimes K}$ on 
      $E^{K}$, $H$ is some function on $E^{K}$ (the Hamiltonian), and 
      $Z_{H}$ is the normalizing constant.
      
      One of the basic example of spin systems is when $E= \{-1,1\}$, 
      and $H= \sum_{A} c_{A}\omega_{A}$, where the functions 
      $\omega_{A}$ are the canonical GKS basis on $\{-1,1\}^{K}$ 
      described before.
      
      To study such systems (and more precisely their asymptotics 
      when $K$ enlarges), one uses some structural inequalities. We 
      present here two fundamental such inequalities, known as GKS 
      inequalities, from Griffiths \cite{Griffiths67}, Kelly and 
      Sherman \cite{KellySherman68}. The GKS property for a basis has 
      been introduced in \cite{bakry-echerbault-GKS}, in an attempt 
      to generalize the GKS inequality to a more general context.
     
     The classical GKS inequalities are settled in the context of 
     $({\Z}/{2\Z})^{K}$.
    As before, we say that $F$ is a GKS function if $F= 
    \sum_{A\subset K} 
    f_{A}\omega_{A}$, where $\forall A$, $f_{A}\geq 0$.
    
    Then we have
    \begin{proposition}
    
    \begin{enumerate} 
    
    \item
    (GKS1 inequality). Assume that $F$ and $H$ are GKS. Then 
    $$\int F\ d\mu_{H} \geq 0.$$
    
    \item (GKS2 inequality). Assume that $F,G$ and $H$ are GKS 
    functions. Then
    $$\int FG\ d\mu_{H} \geq \int F\ d\mu_{H}\int G \ d\mu_{H}.$$
    
    \end{enumerate}
    
    \end{proposition}
    
     The main advantage of the GKS and HGP properties is that 
     they are stable under tensorization. That is, if one considers 
    two sets $(E_{i}, \mu_{i})$ with UOB bases $\mathcal{F}_{i}$ ($i=1,2$), 
    then, on the set 
    $(E_{1}\times E_{2}, \mu_{1}\otimes \mu_{2})$ one has a natural 
    UOB basis $\mathcal{F}_{1}\otimes\mathcal{F}_{2} = (f_{i}\otimes f_{j})$.  Then, if 
    both 
    $\mathcal{F}_{i}$ are GKS or HGP, the same is true for $\mathcal{F}_{1}\otimes 
    \mathcal{F}_{2}$. This is straightforward from the definitions.
    
    This allows us to consider a set $(E, \mu)$ with a given GKS 
    basis $\mathcal{F}$, and then the basis $\mathcal{F}^{\otimes K}$ on $E^{K}$ is 
    again GKS.
    
    One has the following (\cite{bakry-echerbault-GKS})
    
    \begin{proposition}
    If $(E,\mu)$ has a UOB $\mathcal{F}$ which is GKS, then the GKS1 
    inequality is true.
    
    \end{proposition}

    \begin{proof}
    The previous statement just means that if we define a GKS 
    function $F$ as a function which may be written as $F= \sum_{i} 
    F_{i} f_{i}$, where $(f_{i})$ are the elements of $\mathcal{F}$ and 
    $\forall i$, 
    $F_{i}\geq0$, then
    if $F$ and $G$ are GKS functions, one has 
    $$\int F \ d\mu_{H} \geq 0.$$ The statement is straightforward, since
   $\exp (H)$ is again GKS, being the sum of a series with non-negative 
    coefficients, and so $F\exp(H)$ is itself GKS. 
    Since any GKS function has a non-negative integral, the conclusion 
    follows.    \qed
\end{proof}
    
    The GKS2 inequality is much harder. It has only be obtained in 
    some restricted settings, like products of abelian groups, and 
    when the basis comes from    $G^{N}$
    for $G$  elementary groups  like dihedral groups, and some 
  for other few groups. Nevertheless, in any example, one has both 
  the GKS and the HGP property.

  There is no example of a GKS basis where the GKS2 inequality is not 
  satisfied. But we may restrict ourselves to a simpler setting. 
  
  Here is one conjecture that we had been unable to prove, and which 
  motivated most of the material of this section:
  \begin{conjectur}
  If the UOB $\mathcal{F}$ has the GKS and the HGP property, then the GKS2 
  inequality is  true.
  \end{conjectur}

    \section{The hypergroup property in the infinite setting}

 In the first section, we described the hypergroup property in the 
 context of finite setS. In what follows, we consider a general 
 probability space $(E,\mathcal{E}, \mu)$, together with a $L^{2}$ basis 
 $\mathcal{F}= \{1=f_{0}, f_{1}, \ldots, f_{n}, \ldots)$. Very soon, we shall 
 restrict ourselves to the case of a topological space (in fact an 
 interval in the basic examples of Sections \ref{Achour} and 
 \ref{Gasper}), where the functions $f_{i}$ of the basis will be 
 continuous bounded functions. But some general properties may be 
 stated in a more general context.
     
     \subsection{Markov sequences associated with a UOB}
     
     Let $(E,\mathcal{E},\mu)$ a general probability space. In this subsection, we 
     shall ask $(E,\mathcal{E})$ to be at least a ``nice'' measurable space in the 
     context of measure theory.  For us, it shall be enough to suppose 
     that $E$ is a separable complete metric space (a polish space) and 
     that $\mathcal{E}$ is the $\sigma$-algebra of its $\sigma$-field. Then 
     $L^{2}(\mu)$ is separable.

     We suppose that some $L^{2}(\mu)$ orthonormal basis $\mathcal{F}= (f_{0}, 
     \ldots, f_{n}, \ldots)$ is given, with $f_{0}=1$.  In what follows, 
     we shall assume that this is a basis of the real Hilbert space 
     $L^{2}(\mu)$, although we may as well assume that the functions 
     $f_{n}$ may have complex values and be a basis of the complex 
     Hilbert space.  Such a unitary basis $\mathcal{F}$ will be 
     called a Unitary Orthonormal Basis (UOB) associated with the 
     measure $\mu$. 
     
     We are interested in bounded linear operators $K$ on $L^{2}(\mu)$,  
     for which the functions $f_{n}$ are eigenvectors. They are uniquely 
     determined by
     $$K(f_{n})= \lambda_{n} f_{n},$$
     for some bounded sequence $(\lambda_{n})$. The central question we 
     address here is to determine for which sequences $(\lambda_{n})$ one 
     has 
     $$K(f)(x)= \int_{E}f(y)\ k(x,dy),$$
     for some Markov kernel $k(x,dy)$ of probability measures on $E$.
     
   As before, we shall call such a sequence $(\lambda_{n})$ a Markov Sequence  (MS 
     in short) associated with the UOB $\mathcal{F}$.  We shall say 
     that the kernel 
     $k$ (or rather $K$ with a slight abuse of notation) is associated with the MS $(\lambda_{n})$.
     
     In general, this is not an easy question, but as before  the 
     hypergroup property of the basis will be a way of describing all 
     Markov sequences. 
     
     Let us start with some basic remarks.
     
     First, for any Markov operator, $K(f_{0})= f_{0}$, since $f_{0}$ is 
     the constant function, and therefore $\lambda_{0}=1$.
     
     Also, any such Markov operator is symmetric in $L^{2}(\mu)$, since 
     we already know its spectral decomposition which is discrete and 
     given by the basis $\mathcal{F}$. That means that for any pair $(f,g)$ of 
     functions in $L^{2}(\mu)$, one has
     $$\int K(f)(x)g(x)\ \mu(dx) = \int f(x)K(g)(x)\ \mu(dx).$$
     Therefore, the measure $K(x,dy)\mu(dx)$ is symmetric in $(x,y)$.
     
     But any  Markov operator is a contraction in $L^{\infty}(\mu)$, and 
     any symmetric Markov kernel is a contraction in $L^{1}(\mu)$, since, 
     for any $f\in L^{2}(\mu)$,
     $$\int \left|K(f)\right|\ d\mu\leq \int K(\left|f\right|)\ d\mu = \int 
     \left|f\right| K(f_{0})\ d\mu= \int \left|f\right|\ d\mu.$$
     
     Therefore, by interpolation, $K$ is a contraction in $L^{p}(\mu)$ 
     for any $p\geq 1$, and in particular in $L^{2}$.
     
     We deduce from that that any MS $(\lambda_{n})$ satisfies 
     $$\forall n, ~\left|\lambda_{n} \right|\leq 1 .$$
     
     Also, if $(\lambda_{n})$ and $(\mu_{n})$ are MS's, with associated 
     kernels $K$ and $K_{1}$, for any $\theta\in 
     [0,1]$, $(\theta\lambda_{n}+ (1-\theta) \mu_{n})$ is a MS, 
     associated with the kernel $\theta K+ (1-\theta)K_{1}$. 
     
     Therefore, the set of Markov Sequences is  convex, and compact (for 
     the product topology on $\R^{\N}$). This shows that describing all 
     Markov sequences amounts to describe the extremal points of this 
     convex set.
     
     Notice also that the set of all Markov sequences is stable under 
     pointwise multiplication,  which corresponds to the composition of 
     operators. In other words, if $\Theta$ is the set of extremal 
     points in the compact set of Markov sequences, and if 
     $(\lambda_{i}(\theta))$ is the MS associated with the point 
     $\theta\in \Theta$, then one has
     $$\lambda_{i}(\theta)\lambda_{i}(\theta')= \int_{\Theta} 
     \lambda_{i}(\theta_{1})\ R(\theta, \theta', d\theta_{1}),$$ for 
     some probability kernel $R(\theta, \theta', d\theta_{1})$ on the 
     space $\Theta$.

     To determine that an operator $K$ is a Markov operator starting from 
     its spectral decomposition,  we shall need the following proposition.
     
     \begin{proposition} A bounded  symmetric operator $K$ on $L^{2}(\mu)$ is a Markov operator if and only if
     $$K(f_{0})= f_{0}, ~f\geq0 \implies K(f) \geq 0.$$
     
     \end{proposition}

     \begin{proof}
     This is where  we need the fact that the measure 
     space $(E, \mathcal{E}, \mu)$ is a nice space. The conditions on $K$ are obviously 
     necessary.  To see the reverse,  we apply the bi-measure theorem 
     (\cite{dellacherie-meyerT1}, Page 129). We consider the map
     $$\mathcal{E}\times\mathcal{E} \mapsto [0,1] ~: ~(A,B)\mapsto \int\Ind_{A}K(\Ind_B)\ 
     d\mu(x).$$
      For any fixed $B$, this is a measure in $A$, and by symmetry, it is 
      also a measure in $B$. Since our spaces are polish spaces, theses 
      measures are tight, and therefore we may extend this operation into 
      a measure $\mu_{K}(dx,dy)$ on the $\sigma$-algebra $\mathcal{E}\times 
      \mathcal{E}$.   The measure is symmetric, and any of its marginal is 
      $\mu$.
      
      Then, we apply the measure decomposition theorem to write 
      $$\mu_{K}(dx,dy)= K(x,dy) \mu(dy).$$
      The kernel $K(x,dy)$ is exactly the kernel we are looking for.\qed
\end{proof}
\subsection{The hypergroup property}

The GKS property is relatively easy to state in a general context, as 
soon as the functions $f_{i}$ of the basis are in $L^{3}(\mu)$, since 
then we may just ask that
$$\forall i,j,k, ~ \int_{E} f_{i}f_{i}f_{k}\ d\mu \geq 0.$$

But for the dual hypergroup property, one has to be a bit more 
cautious. In general, functions in $L^{2}$ are defined up to a set of 
$\mu$-measure $0$. Therefore, the meaning of 
${f_{i}(x)}/{f_{i}(x_{0})}$ is not so clear. In order to avoid 
difficulties, and since this shall correspond to the examples we are 
going to describe below, we restrict ourselves to the following 
setting: $E$ is a compact separable Hausdorff space, and the 
functions $f_{i}$ are continuous on $E$.

We may then set the following definition.

\begin{definition}
     We shall say that the UOB $\mathcal{F}$ has the hypergroup property (HGP in short) at 
     some point $x_{0}\in E$, if, for any $x\in E$, the operator 
     defined on $\mathcal{F}$ by
     $$K_{x}(f_{i}) = \frac{f_{i}(x)}{f_{i}(x_{0})} f_{i}$$ is a 
     Markov operator.
     \end{definition}
     
     In other words, we require the sequences 
     $\big({f_{i}(x)/}{f_{i}(x_{0})}\big)_i$ to be Markov sequences.
     
     Observe that this implies that $\left|f_{i}\right|$ is maximal 
     at $x_{0}$, since the eigenvalues of a Markov operator must be 
     bounded by 1. Since the functions $f_{i}$ are normalized in 
     $L^{2}(\mu)$, then for any $i\in \N$, one has  $\left|f_{i}(x_{0})\right|\geq 1$.

     As before, this definition is equivalent to the following
     
     \begin{proposition}\label{HGP2}
     
    The UOB $\mathcal{F}$ has the HGP property at the point $x_{0}$ if and only if there exists a 
    probability kernel $K(x,y,dz)$ such that, for any $i\in \N$
    $$\frac{f_{i}(x)f_{i}(y)}{f_{i}(x_{0}) }= \int 
    f_{i}(z)\ K(x,y,dz).$$
     \end{proposition}

     \begin{proof}
     
     We shall mainly use this in the obvious way: if there is a 
     probability kernel $K(x,y,dz)$ satisfying the hypothesis of the 
     proposition, then the HGP property holds. 
     In fact, if such a probability kernel $K(x,y,dz)$ exists, for any 
     $x\in E$, the Markov kernel $K(x,y,dz)$ defines a Markov 
     operator with Markov sequence $\big({f_{i}(x)}/{f_{i}(x_{0})}\big)_i$.
     
     For the reverse, if the HGP property holds, there exists for any 
     $x$ a Markov kernel $k_{x}(y,dz)$ which satisfies 
     $$\frac{f_{i}(x)f_{i}(y)}{f_{i}(x_{0}) }= \int 
     f_{i}(z)\ k_{x}(y,dz).$$
     
     It remains to turn this  family of kernels into a two parameters 
     kernel $K(x,y,dz)$.\qed
\end{proof}

     Sometimes, it is easier to see the HGP property in the reverse 
     way. 
     
     \begin{proposition}\label{HGP3} 
     If there exists a non-negative kernel $k(x,dy,dz)$ such that
      for any $i,j$
     $$\int f_{i}(y)f_{j}(z) \ k(x,dy,dz) 
   = \delta_{ij} \frac{f_{i}(x)}{f_{i}(x_{0})}$$
     then the HGP property holds at the point $x_{0}$.
     
     \end{proposition}
     
     \begin{proof} Remark first that from our hypotheses,  the kernel 
     $k(x,dy,dz)$ is a probability kernel (taking $i=j=0$ in the 
     definition).
     
     Now, the two marginals of the kernel $k(x,dy,dz)$ are equal to 
     $\mu$, since
     $$\int f_{i}(y)\ k(x,dy,dz)= \delta_{0i},$$ which shows that 
     those marginals and $\mu$ give the same integral to any $f_{i}$, 
     and therefore to any $L^{2}$ function.

     We may then decompose the kernel $k(x,dy,dz)= 
     k_{1}(x,y,dz)\mu(dy)$.
     
     Since the functions $f_{i}$ are bounded, we may consider the 
     bounded functions 
     $$H_{i}(x,y)=\int f_{i}(z)\ k_{1}(x,y,dz).$$
     From the definition of $k_{1}$ and the hypothesis on $k$, it is 
     straightforward to check that
     $$\int H_{i}(x,y) f_{k}(x)f_{l}(y)\ \mu(dx)\mu(dy)= 
     \delta_{il}\delta_{ik} \frac{1}{f_{i}(x_{0})},$$ and hence
     $$H_{i}(x,y) = \frac{f_{i}(x)f_{i}(y)}{f_{i}(x_{0})}.$$
     
     Therefore $k_{1}$ satisfies the hypotheses of Proposition 
     \ref{HGP2} and the proof is completed.\qed
\end{proof}
     
     The representation of Markov sequences is then the same than in 
     the previous section.
     
     \begin{theorem} 
     
     If the UOB $\mathcal{F}$ has the hypergroup property, then 
    any Markov sequence has the representation
    $$\lambda_{i}= \int _{E} \frac{f_{i}(x)}{f_{i}(x_{0})}\  
    \nu(dx),$$ for some probability measure $\nu$ on $E$. 
    Moreover, the Markov sequences $$\left(\frac{f_{i}(x)}{f_{i}(x_{0})}\right)_i$$ 
    are the extremal Markov sequences.
     
     \end{theorem}
     
     \begin{proof} The proof is exactly similar to the finite case (see Theorem 
     \ref{thmfond}).\qed
\end{proof}

Remark that the series $$K(x,y,z)=\sum_{i} 
\frac{f_{i}(x)f_{i}(y)f_{i}(z)}{f_{i}(x_{0})}$$ does not converge in 
general. We shall see in the examples developed in the next section that the formal measure 
$K(x,y,z)\mu(dz)$ is not absolutely continuous with respect to $\mu$, 
and may have Dirac masses at some points.

But we may still define a convolution structure from $f_{i}*f_{j}= 
\delta_{ij}{f_{i}}/{f_{i}(x_{0})},$ which maps probability 
measures onto probability measures, and all Markov kernels associated 
with $\mathcal{F}$ would 
be represented as $K(f) = f* \nu$, for some probability measure on 
$E$. We give no details here since this will not be used in the sequel.

     \section{Sturm--Liouville bases and Achour--Trim\`eche's theorem}
     
    \subsection{The natural UOB associated with a measure on a 
    compact interval}
    
    In this section, we shall consider some  natural infinite UOB 
    coming from the spectral decomposition of  Sturm--Liouville 
    operators on a compact interval of the real line.
    
    Let us first describe the context. Consider a probability measure $\mu(dx)= 
    \rho(x)dx$ on some compact interval $[a,b]\subset \R$. In what 
    follows, we shall assume for simplicity that $\rho$ is smooth, 
    bounded above and away from $0$ on $[a,b]$. The density 
    $\rho$ is associated with a canonical differential operator 
    $$L(f)(x) = f''(x) + \frac{\rho'}{\rho}(x) f'(x),$$
    which is symmetric in $L^{2}(\mu)$. We shall consider here $L$ 
    acting on functions on $[a,b]$ with derivative $0$ at the 
    boundaries $a$ and $b$ (Neumann boundary conditions).
    
    In this context, $L$ is essentially self adjoint on the space of 
    smooth functions with $f'(a)=f'(b)=0$ and  there is an 
    orthonormal basis 
    $$\mathcal{F}= (1=f_{0}, f_{1}, \ldots, f_{n},\ldots)$$ of $L^{2}(\mu)$ 
    which is given by eigenvectors of $L$ satisfying the boundary 
    conditions. This means that there is an increasing sequence of 
    real numbers $$0=\lambda_{0}<\lambda_{1}<\cdots< 
    \lambda_{n}<\cdots$$ such that
    $$Lf_{i}= -\lambda_{i} f_{i},~f_{i}'(a)= f_{i}'(b)=0.$$
    
    From the standard theory of Sturm--Liouville operators, the 
    eigenvalues $\lambda_{i}$ are non-negative and simple. Therefore, 
    there is for any $\lambda_{i}$ a unique solution $f_{i}$ of the 
    previous equation, which has norm 1 in $L^{2}(\mu)$ and which 
    satisfies $f(a) > 0$. We refer to any standard text book for 
    details (see \cite{brezis} or \cite{zettl05} for example).
    
    This basis shall be called the canonical UOB associated with 
    $\mu$ on $[a,b]$.
    
    The fact that we chose to deal with the Numen boundary 
    conditions and not with the Dirichlet boundary conditions 
    ($f(a)=f(b)=0$) comes from the fact that we require the function 
    $1$ to be an eigenvector of the operator.

    It will be much more convenient in what follows, essentially for 
    notations, to extend our functions by symmetry in $a$ 
    and $b$, (and the same for $\mu$). In this way we may consider 
    that we are working on functions on the real line, which are 
    symmetric under $x\mapsto 2a-x$, and are $2(b-a)$-periodic.
    
    The eigenvectors are perhaps not smooth then at the boundaries 
    $a$ and $b$, but they are at least $C^{2}$ (since they are 
    solutions of the equation $Lf_{i}= -\lambda_{i} f_{i}$ at the 
    boundaries).

    The hypergroup property is stated at some point in $[a,b]$. In 
    the finite case, we know at which point we may expect the 
    hypergroup property to hold: this is a point of minimal mass. In 
    the general case, such a reasoning does not hold, since one may 
    choose a point with minimal density to a given reference measure, 
    but this depends on this choice.  
    
    Here, the basis $\mathcal{F}$ is the 
    sequence of eigenvectors of an elliptic second order 
    differential operator $L$ symmetric in $L^{2}(\mu)$. In this setting, there 
    is a natural distance associated with the operator $L$ (in this 
    precise example of Sturm--Liouville operators, this is the natural 
    distance on $\R$). In any example we know, the point $x_{0}$ is 
    minimal in the following sense
    \begin{equation}\label{massmin}\lim_{r\to 0} 
    \frac{\mu\big(B(x_{0},r)\big)}{\mu\big(B(x,r)\big)} \leq 1.\end{equation}
    
    We did not try to prove this in a more general context. However, 
    it is not clear how the properties of the operator $L$ must be 
    reflected in the properties of $\mathcal{F}$ to insure for example that 
    the maximal values of the eigenvectors are attained at the same 
    point, and that this point is of minimal mass in the sense of 
    \eqref{massmin}.

     \subsection{Wave equations}
     In this context, one has some other interpretation of the 
     hypergroup property.
     
     On $D=[a,b]^{2}$, we shall consider the following differential 
     equation
     \begin{equation} \label{WaveEq} L_{x}F(x,y)= L_{y}F(x,y),\end{equation} for a function $F$ which has 
     Neumann boundary conditions on the boundary of $D$. We shall say 
     that such a function is a solution of the (modified) wave 
     equation. 
     
     We have to be careful here with the regularity assumption on the 
     function $F(x,y)$ that we require. We shall see later that given any 
     smooth function $f(x)$ at the level $x=x_{0}$, with 
     Neumann boundary conditions, there is exactly one smooth function 
     $F(x,y)$ on $D$ 
     which is solution of  Equation \eqref{WaveEq} and satisfies 
     $F(x,y_{0})= f(x)$. 
     
     In fact, if $f(x)= \sum_{i} a_{i} f_{i}(x)$ is the $L^{2}$ 
     orthogonal decomposition of $f$, then 
     $$F(x,y)= \sum_{i} \frac{a_{i} }{f_{i}(y_{0})} f_{i}(x)f_{i}(y)$$ 
     is a formal $L^{2}$ solution of the wave 
     equation,  since 
     $$L_{x}F= \sum_{i} \lambda_{i}\frac{a_{i} 
     }{f_{i}(y_{0})} f_{i}(x)f_{i}(y)= L_{y} F.$$
     But we do not even know (for the moment) that this solution is 
     such that $L_{x} F$ is in $L^{2}(\mu\otimes \mu)$.
     
     Therefore, we shall say that $F$ is a weak $L^{2}$ solution of 
     \eqref{WaveEq} if for any smooth function $G(x,y)$ with Neumann 
     boundary conditions on $\partial D$, one has
     $$ \int [(L_{x}-L_{y}) G(x,y)] F(x,y)\ \mu(dx)\mu(dy) =0.$$
     
     Since $$\int L_{x} (F) G\ \mu(dx)\mu(dy) = \int L_{x} (G) F \ 
     \mu(dx) \mu(dy) $$ for any pair of smooth functions $F$ and $G$ 
     satisfying the Neumann boundary conditions, then any ordinary 
     solution is a weak one.
     
     Now, given any $L^{2}(\mu)$ function $f(x)$, the above 
     construction produces a weak $L^{2}(\mu\otimes \mu)$ solution 
     $F(x,y)$ satisfying the wave equation \eqref{WaveEq}, and we 
     claim  immediately that this solution is unique. In fact, writing 
     the function $F(x,y)$ as 
     $$F(x,y) = \sum_{ij} a_{ij} f_{i}(x) f_{j}(x),$$
    and using the fact that the eigenvalues are simple, one may check 
    that if $F$ satisfies weakly  \eqref{WaveEq}, then $a_{ij}= 0$ if $i\neq 
    j$, from which we deduce  our claim. 
    
    Observe moreover that $a_{00}= \int F(x,y)\ \mu(dx)\mu(dy)$, and 
    that for almost every $y_{0}$, $F(x,y_{0})\in L^{2}(\mu)$ and that
    $$\int F(x,y_{0})\ \mu(dx)= \int F(x,y)\ \mu(dx) \mu(dy).$$
    
    It is not clear however that if $F(x, y_{0})$ is smooth, then 
    $F(x,y)$ is smooth. This shall be done later at least when 
    $y_{0}=a$ or $y_{0}=b$.
    
    The link between solutions of the wave equation and Markov kernels 
    is the following. 
    
    If a Markov kernel is Hilbert--Schmidt (that is if its 
    eigenvalues $\lambda_{i}$ satisfy  $\sum_{i} \lambda_{i}^{2} < 
    \infty$), then it may be represented as a 
    $$K(f)(x)= \int f(y) k(x,y)\ \mu(dy),$$
     where 
     $$k(x,y)= \sum_{i} \lambda_{i} f_{i}(x) f_{i}(y).$$
     Therefore, there is a one-to-one correspondence between 
     Hilbert--Schmidt Markov kernels and non-negative  weak $L^{2}$ solutions of 
     the wave equation which satisfy
     $$\int F(x,y)\ \mu(dx)\mu(dy)= 1.$$
     
     We then have the following
     
     \begin{theorem}\label{ThmWE1} Assume that for any function $f(x)$ on the 
     interval $[a,b]$ 
     with Numen boundary conditions, there exists a unique $C^{2}$ 
     solution $H_{f}(x,y)$ 
     of the wave equation \eqref{WaveEq} on $[a,b]^{2}$ such that 
     $H(x,y_{0})= f(x)$.  
     Then,  the HGP 
     property holds at the point $y_{0}$ for the natural UOB 
     associated with $\mu$ if and only if
     whenever $f\geq0$ one has $H_{f}\geq0$ on $I^{2}$.
    \end{theorem}

    In other words, the hypergroup property is equivalent to the fact 
    that the wave equation is positivity preserving.
 
 \begin{proof} Assume first that the hypergroup property holds at the point 
 $y_{0}$.  Take any smooth solution $F(x,y)$ of the wave equation 
 with $F(x,y_{0})= f(x)\geq 0$. Then, from what we just saw, one has
 $$F(x,y) = K_{y}(f)(x),$$
  where $K_{y}$ is the Markov kernel with eigenvalues 
  ${f_{i}(y)}/{f_{i}(y_{0}) }$. Therefore, $F(x,y)$ is everywhere 
  non-negative.
  
  On the other hand, assume that any smooth solution of the wave 
  equation which is non-negative on $\{y=y_{0}\}$ is non-negative 
  everywhere. Consider the heat kernel 
  $$p_{t}(x,z) = \sum_{i} \exp(-\lambda_it)f_{i}(x)f_{i}(z).$$
  We know that it is a smooth function on $D$, which is everywhere 
  positive. Then,
  $$F_{t,z}(x,y)= \sum_{i} 
  \frac{\exp(-\lambda_{i}t)}{f_{i}(y_{0})}f_{i}(z)f_{i}(x)f_{i}(y)$$ 
  is the unique $L^{2}$  
  solution of the wave equation with $F_{t,z}(x,y_{0}) = p_{t}(x,z)$.

  Therefore, this function is non-negative, and this shows that for 
  any $t>0$, the sequence 
  $$\exp(-\lambda_{i}t) \frac{f_{n}(z)}{f_{n}(y_{0})}$$ is a Markov 
  sequence. It remains to let $t$ go to $0$ to get the result, since 
  a limit of Markov sequences is a Markov sequence.\qed
\end{proof}
  
  \subsection{Achour--Trim\`eche's theorem and wave equations}\label{Achour}
  
  In what follows, we consider the case  of a symmetric interval 
  $[-b,b]$. Then we have
  \begin{theorem}[Achour--Trim\`eche]
  Let  $\rho$ be a log-concave and symmetric density on 
  $[-b,b]$. Then, the natural UOB associated with $\mu$ has the HGP 
  property at the point $-b$.  In this  case, we may as well choose 
  $x_{0}=b$.
  The same is true on any interval with any log-concave increasing 
  density $\rho$. 
  
  \end{theorem}
  
  This result is one of the very few cases when one may produce 
  hypergroup bases without any kind of group structure on the space 
  $E$. We shall see in the next chapter that this property holds for 
  Jacobi polynomials, but in this case, there are at least for the 
  integer values of the coefficients some interpretations of the 
  convolution which reflects the group action of some orthogonal 
  group. There is absolutely no such interpretation in this context.
  
  In general, Achour--Trim\`eche's result is stated with a density $\rho$ which 
  vanishes on the boundary. Under the conditions usually stated in 
  Achour--Trim\`eche's theorem, there are then no difference between Neumann and 
  Dirichlet boundary conditions.  The series $$\sum_{i} 
  \frac{f_{i}(x)f_{i}(y)f_{i}(z)}{f_{i}(x_{0})}\mu(dz)$$  is absolutely 
  continuous with respect to the measure $\mu$, which is not the case 
  here.
  
  We chose to present this result in the case where the density 
  $\rho$ is bounded from below because it seemed to us to be more natural.
  
  Apparently, the proof of Achour--Trim\`eche's theorem had never been published. We  
  found a mention of it in the reference book \cite{BloomHeyer95} and the result is announced in \cite{AchourTrimeche}, with 
  no proof. Most of the ideas presented here come from Achour's thesis.  The idea follows 
  a previous result of Chebli \cite{Chebli}, which works on $[0, \infty)$ and is somehow simpler (It corresponds to the case of a concave decreasing density). 
  
  \begin{proof}
  To prove this result, we shall make use of the characterization of 
  the hypergroup property in terms of the wave equation given in 
  Theorem  \ref{ThmWE1}.
  We shall see in the next paragraph that any smooth bounded function 
  satisfying Neumann boundary conditions on $[a,b]$ has a unique extension 
  as a smooth solution of the  wave equation \eqref{WaveEq}, when 
  $x_{0}$ is one of the boundary points (see Paragraph \ref{more}). 
  (This has nothing to do with the log-concavity of the measure or 
  with the symmetry: this is just a consequence of the fact that 
  $\log (\rho)$ is smooth and bounded.) 
  
  We first treat the case where the density $\rho$ is log-concave 
  symmetric.
   
  First we make use of the symmetry assumption. Then, any eigenvector 
  of the operator $L$ on $[-b,b]$ with Neumann boundary conditions is 
  either even or odd, since $f_{i}(-x)$ is also an eigenvector with 
  the same eigenvalue.
  
  Then, any $L^{2}$ solution $F(x,y)$ of the wave equation 
  \eqref{WaveEq}, written as
  $$F(x,y)= \sum_{i} a_{i}f_{i}(x)f_{i}(y)$$is symmetric under the 
  change $(x,y)\mapsto (y,x)$ and under $(x,y)\mapsto (-x,-y)$.
  
  We want to show that if $F(x,-b) \geq 0$, then $F(x,y)\geq 0$ 
  everywhere. For this, it is enough to show this on the domain 
  $D_{1}=\{x+y\leq 0, ~x\geq y\}$. 
  
  Also, we may change $F$ into $F+ epsilon$ for any $epsilon>0$, and 
  we are thus reduced to prove that the result is true when the 
  function $f$ on the boundary is bounded below by some positive 
  constant.

  \begin{figure}[ht]
  \centering
		\includegraphics[width=.8\linewidth]{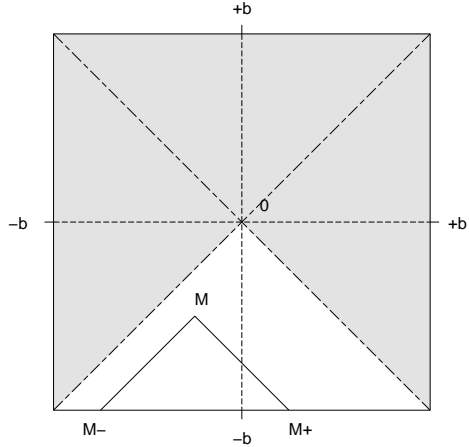}
		\caption{Triangle $\Delta_M$.}
		\label{fig:MM+M-}
\end{figure}
  
  Then, a point $M\in D_{1}$, let $\Delta_{M}$ be the triangle 
  delimited by the lines $x+y=c$ and $x-y=c'$ passing through $M$ and 
  the line $y= -c$.  Let $M_{-}$ and $M_{+}$ be the points of this 
  triangle which lie on the line $\{y=-b\}$, $M_{-}$ being the left 
  point and $M_{+}$ the right one (see Figure \ref{fig:MM+M-}).  Let $F$ be a smooth 
  solution of the wave equation \eqref{WaveEq}, and let $G(x,y)= 
  F(x,y)\rho(x)\rho(y)$.
  
  \begin{equation}\label{repr1}
  2G(M)= G(M_{-})+ G(M_{+}) + \int_{[M_{-}M]} G(s) a_{+}(s)\ ds + 
  \int_{[M M_{+}]} G(s) a_{-}(s)\ ds
  \end{equation}
  where 
  $$a_{+}(x,y) = \frac{1}{\sqrt{2}}\left(\frac{\rho'}{\rho}(x)+ 
  \frac{\rho'}{\rho}(y)\right), 
  ~a_{-}(x,y)= 
  \frac{1}{\sqrt{2}}\left(\frac{\rho'}{\rho}(y)-\frac{\rho'}{\rho}(x)\right),$$
   and the integral $\int_{[M_{-}M]}H(s)\ ds$ and 
   $\int_{[MM_{+}]}H(s)\ ds$ denote the one dimensional integrals along 
   the segments $[M_{-}M]$ and $[MM_{+}]$ against the (euclidean) length
   measure on those lines.
   
   This formula relies on an integration by parts formula. Even 
   though we shall use it only on domains like $\Delta_{M}$, it is 
   perhaps of some interest to state it in general. So we set it as a 
   lemma. We shall not give too much details here, since a more 
   general formula will be derived in the next paragraph.
   
   \begin{lemma} Let $H$ a smooth function on $D$ and 
    $\Omega\subset D$ a domain  with a 
   piece-wise $C^{1}$ boundary  $\partial \Omega$. Then
   $$\int_{\Omega}(L_{x}-L_{y})H(x,y) ~\rho(x)\rho(y)\ dx dy= 
   \int_{\partial\Omega} \nabla H \odot n ~\rho(x)\rho(y)\ ds,$$
 where $n=(n_{x},n_{y})$ denotes the exterior normal derivative of 
 the domain and 
 $$\nabla H\odot n= \partial_{x} H n_{x}- \partial_{y} H n_{y},$$
 $ds$ designing the length measure on the boundary $\partial\Omega$.
 \end{lemma}
 
 We shall not prove this lemma. It is the analogue of the classical 
 Stokes formula, where the elliptic operator $\Delta$ is replaced by 
 the hyperbolic operator $L_{x}-L_{y}$, its invariant measure being 
 $\rho(x)\rho(y)dxdy$.

 One may see this as a particular case of the general integration by 
 parts formula
 $$\int_{D} H_{1}(L_{x}-L_{y})H_{1}~\rho(x)\rho(y)\ dxdy= -\int \nabla 
 H_{1}\odot \nabla H_{2}~ \rho(x)\rho(y)\ dx dy,$$ applied with $H_{1}= 
 \Ind_{\Delta_{M}}$.

 From the previous formula,  applied on $\Omega = \Delta_{M}$ for a function $F$ which is solution of 
 the wave equation \eqref{WaveEq} and has normal derivative vanishing on the 
 boundary $[M_{-},M_{+}]$, one has
 $$-\int_{[M-M]} (\partial_{x}F+ \partial_{y}F)\rho(x)\rho(y)\ ds + 
 \int_{[MM_{+}]} (\partial _{x} F - \partial_{y} F) \rho(x)\rho(y)\ ds 
 =0.$$

 We may then perform a next integration by parts on both integrals 
 to find
 \begin{multline*}
G(M_{-})-G(M) + \int_{[M_{-}M]} G(s) a_{+}(s)\ ds \\+ G(M_{+})-G(M) + 
 \int _{[MM_{+}]}G(s) a_{-}(s)\ ds =0,
\end{multline*}
 which gives \eqref{repr1}.
 
 Under our assumptions of $\rho$, both $a_{+}$ and $a_{-}$ are non-negative on the subdomain $D_{1}$:  under the 
 log-concavity assumption $a_{-}$ is non-negative on $\{y\leq x\}$, 
 and $a(x)+a(y)= a(y)-a(-x) \geq 0 $ if $x+y \leq0$.
 
 Now, consider the smallest $y\in (-b,0)$  such that there exists 
 some point in $(x,y)\in D_{1}$ with $G(x,y) =0$.  On this point, we 
 have 
 $$2G(M)=0 \geq G(M_{-})+ G(M_{+}),$$
 which gives a contradiction.

 For the case where the density is log-concave increasing,  we 
 may use the same argument on the domain $\{x\geq y\}$, since we still have the solution of the wave equation symmetric under the change $(x,y)\mapsto (y,x)$. Then we extend $F$ by symmetry around the axes $x= -b$ and 
 $x=b$, and then by periodicity, into a function defined on $\R\times 
 [-b,b]$. The same argument of integration by parts remains valid, and, 
 by means of the symmetrization, the domain of integration   $(Mm_+M_-)$ that should be used  is replaced by  the same triangle $\Delta_M$ as before, as  shown in Figure \ref{fig:MM+M-pos}. Then we use the fact  that the function $a$ is decreasing and  non-negative.
 
   \begin{figure}[ht]
   		\centering
		\includegraphics[width=.8\linewidth]{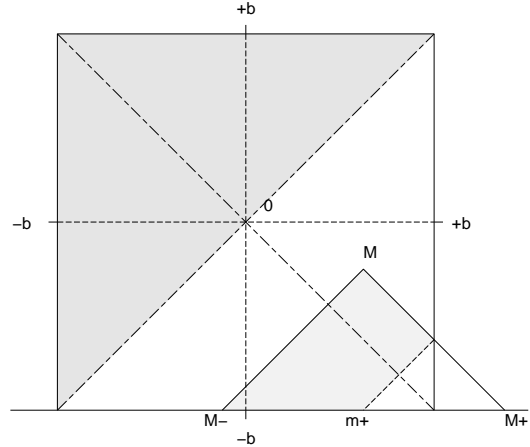}
		\caption{Triangle $\Delta_M$ and $(Mm_+M_-)$.}
		\label{fig:MM+M-pos}
\end{figure}\qed
\end{proof}
 
 Notice that the second case (when $a$ is non-negative) may be reduced to the first one if we 
 extend $\rho$ by symmetry around $b$, into a log-concave  
 function on the interval $[-b,3b]$, symmetric around the point $x=b$.
 (The function $\rho$ may not be $C^{2}$ at the point $x=b$, but this 
 causes no problem).  Then one has to apply the previous  result on symmetric 
 functions on the interval $[-b,3b]$.

 \subsection{Other representations of the solutions of the wave 
 equation}
 
 In this paragraph, we shall consider an operator $L(f)= f''+ a(x) 
 f'$ in $I=[0,1]$, with Neumann boundary conditions, and look at 
 different representations of the solutions of the wave equation
\eqref{WaveEq}.
 
 We shall consider the probability measure $\rho(x)dx$ in $I$ which 
 satisfies ${\rho'}/{\rho}= a$.
 
Let  $F(x,y)$ be a solution of the wave equation $(L_{x}-L_{y})F=0$ 
on $I^{2}$ with Neumann boundary conditions.  As before, it is easier to 
 extend   $F$  to 
 $\R^{2}$ by imposing symmetry conditions at the boundaries $x\in \Z$ 
 or $y\in \Z$, and to extend $a$ by imposing antisymmetry conditions 
 on these lines (or if one prefers, symmetry conditions on $\rho$).
 
 We have seen before that such an equation has an integral 
 representation \eqref{repr1}. Our first task shall be to change it 
 into a new one. 
 
 As before, for $M=(X,Y)$ with $y>0$, we denote by  $\Delta_{M}$ the 
 triangle delimited by the lines $\{x+y= X+Y\}, \{y-x= Y-X\}$ and 
 $\{y=0\}$. $M_{+}$ and $M_{-}$ denote the edges of this triangle 
 which lie on the line $\{y=0\}$.
 
\begin{figure}[ht]
   		\centering
		 \includegraphics[width=.8\linewidth]{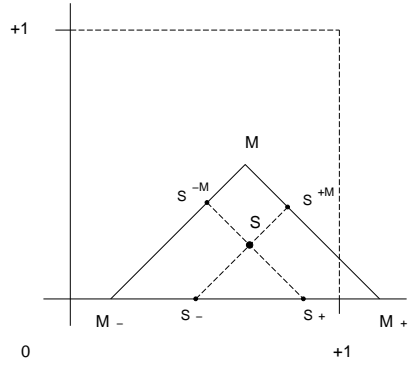}
		 \caption{$S$, $S_+$, $S_-$, $S^{+M}$ et $S^{-M}$.}
		 \label{fig2}
 \end{figure}

 For $S\in \Delta_{M}$,  we denote by $S^{-M}$ the unique 
 point $U$ on the interval $[M_{-},M]$ such that $S\in [U,U_{+}]$ and $S^{+M}$ the 
 unique point $U$ on $[M ,M_{+}]$ such that $S\in [U_{-}U]$ (see  Figure \ref{fig2}). 

 Recall 
 that 
 $$a_{+}(x,y)= \frac{1}{\sqrt{2}}\big(a(x)+a(y)\big), ~a_{-}(x,y)= 
 \frac{1}{\sqrt{2}}\big(a(y)-a(x)\big),$$
 and  let $R(x,y)=\rho(x)\rho(y)$.
 
 \begin{proposition} \label{repr2} If a continuous function $G$ satisfies \eqref{repr1}, then the 
 function $$H(x,y) = \frac{G(x,y)}{\sqrt{R(x,y)}}$$ satisfies
 
 \begin{multline*}
2H(M)= H(M_{-})+ H(M_{+}) \\+ \int_{[M_{-},M_{+}]}H(S) a_{0}(M,S)\ dS+ 
 \int_{\Delta_{M}}  H(S)a(M,S)\ dS,
\end{multline*}
  where
  $$a(M,S)= \frac{1}{2}\sqrt{R(S)}\left(\frac{a_{+}(S^{-M})}{\sqrt{R(S^{-M})}}a_{-}(S)+ 
  \frac{a_{-}(S^{+M})}{\sqrt{R(S^{+M})}}a_{+}(S)\right),$$ and 
  $$a_{0}(M,S)= 
  \frac{1}{2\sqrt{2}} \sqrt{R(S)}
  \left(\frac{a_{+}(S^{-M})}{\sqrt{R(S^{-M})}}+\frac{a_{-}(S^{+M})}{\sqrt{R(S^{+M})}}\right).$$
  \end{proposition}
 \begin{proof}
 We first start by considering the function 
 $$\theta_{M}(S)= \sqrt{\frac{R(M)}{R(S)}},$$
 and we notice that, for $S\in [M_{-},M]$,
 $$\theta_{M}(S)= 1+ \frac{1}{2}\int_{[S,M]} \theta_{M}(U) a_{+}(U)\ 
 dU,$$ and similarly that, for $S\in[M,M_{+}]$,
 $$\theta_{M}(S)= 1+ \frac{1}{2}\int_{[M,S]} \theta_{M}(U) a_{-}(U)\ 
  dU.$$
  
  Then starting from Equation \eqref{repr1}, we  replace the term
  $$\int_{[M_{-}M]}G(S)a_{+}(S)\ dS$$ by
  $$\int_{[M_{-}M]}(1-\theta_{M}(S))G(S)a_{+}(S)\ dS+ 
  \int_{[M_{-}M]}\theta_{M}(S)G(S)a_{+}(S)\ dS.$$
  In the last integral, replace 
  $G(S)$ by 
  $$\frac{1}{2}\left(G(M_{-})+ G(S_{+})+ \int_{[M_{-}S]}G(U)a_{+}(U)\ dU+ 
  \int_{[SS_{+}]} G(U)a_{-}(U)\ dU\right).$$
  Then we have 
  \begin{eqnarray*}
\lefteqn{
\int_{S\in [M_{-}M]}\theta_{M}(S) a_{+}(S)
  \left(  \int_{U\in[M_{-}S]}  G(U) a_{+}(U)\ dU  \right)\
  dS =
  }&\phantom{on decale}&
   \\
  &&
  \int_{[M_{-}M]}G(S)a_{+}(S) 
  \left(\int_{U\in [S,M]}  \theta_{M}(U)a_{+}(U)\ dU\right)\ 
  dS,
\end{eqnarray*}
  
  and this last expression cancels with
  $$\int_{[M_{-}M]}(1-\theta_{M}(S)) G(S) a_{+}(S)\ dS.$$
  
  We do the same computation on the other side $[M,M_{+}]$, and we 
  collect the results. Observe that the term
  $$\int_{[M_{-}M]} \theta_{M}(S) a_{+}(S) G(S_{-})\ dS$$ gives rise 
  to one part of the integral
  $$\int_{[M_{-}M_{+}]} H(U) a_{0}(M,U)\ dU,$$
   while the term
   $$\int_{[M_{-}M]} \int_{[SS_{+}]} G(U) \theta_{M}(S) a_{+}(S) 
   a_{-}(U)\ dSdU$$ produces one part of the integral
   $$\int_{\Delta_{M}} H(U) a(M,U)\ dU.$$\qed
\end{proof}
   
   We shall see in what follows that one may find many other integral 
   representations of the wave equation.
   
   One is the following
   
   \begin{proposition}\label{proppsi} Let $F$ be a solution of the wave equation \eqref{WaveEq}, and 
   let $\psi$ be a smooth positive function on $I^{2}$ satisfying the Neumann 
   boundary conditions.
   Then, if we set $G(x,y)= F(x,y)\psi(x,y) \rho(x)\rho(y)$, we have
   \begin{multline*} 2G(M)=  G(M_{+})+ G(M_{-})+ \int_{[M_{-}M]}G(s) K_{+}(s)\ ds \\
   +   \int_{[MM_{+}]}G(s) K_{-}(s)\ ds - \int_{\Delta_{M}} 
   G(s)\frac{(L_{x}-L_{y})(\psi)}{\psi}\ ds,
   \end{multline*}
  where $$K_{+}(x,y)= \sqrt{2}(\partial_{x}+\partial_{y}) 
  \log\left(\psi \sqrt{\rho(x)\rho(y)}\right)$$ and 
  $$K_{-}(x,y)= \sqrt{2}(\partial_{y}-\partial_{x}) 
  \log\left(\psi \sqrt{\rho(x)\rho(y)}\right).$$
   \end{proposition}

   \begin{proof} The proof is the same as before, but we have to consider the 
   equation satisfied by $G$ instead of the equation satisfied by 
   $F$.  It is perhaps easier to make the computations under a change 
   of variables 
   $$x=\frac{u+v}{\sqrt{2}}, ~y= \frac{u-v}{\sqrt{2}},$$ in which 
   case the operator $L_{x}-L_{y}$ becomes 
   $$2\partial^{2}_{uv}-\frac{a_{-}(u,v)}{\sqrt{2}}\partial_{u}+\frac{a_{+}(u,v)}{\sqrt{2}}\partial_{v}.$$

   We extend our functions by symmetry to the set $\{y\leq 0\}$, and 
   then our functions become for the new variables symmetric under the symmetry 
   $(u,v)\mapsto (v,u)$.
   Then the result is obtained through the integration on a square 
   $\{u_{0}\leq u\leq u_{1}, ~u_{0}\leq v \leq u_{1}\}$. Since we 
   shall not use these representations here, the details are left to the 
   reader.\qed
\end{proof}
   
   As a consequence, if we set $U\leq V$ when $U\in \Delta_{V}$, and 
     if there exists a positive function $\psi$ satisfying the Neumann 
     boundary conditions with $(L_{x}-L_{y})(\psi)\leq 0$ and such that 
     $\sqrt{\rho(x)\rho(y)}\,\psi$ is increasing for this partial order, then any continuous 
     solution of the wave equation \eqref{WaveEq} which is non-negative on 
     $\{y=0\}$ is non-negative everywhere. In particular, if there is a 
     solution of the wave equation which is increasing for this order, 
     then the property holds.

 If we are looking for the hypergroup property  at the point $0$ for the Neumann basis 
 on $[0,1]$, a good candidate for the function $\psi$ in Proposition 
 \ref{proppsi} seems to be 
 $$\psi(x,y)= 1-\frac{f(x)f(y)}{f(1)^{2}},$$
 where $f$ is the (increasing) eigenvector associated with the first 
 non-$0$ eigenvalue, provided that $\left|f(1)\right|\leq 
 \left|f(0)\right|$, which is a necessary condition for the 
 hypergroup property to hold at $0$.
 But we were unable to derive reasonable conditions on $a$ which 
 would insure that for this particular case the function $\sqrt{\rho(x)\rho(y)}\psi$ is increasing for the 
 partial order on $[0,1]^{2}$.

 \subsection{More about the solutions of the wave equation 
 \eqref{WaveEq}}\label{more}
 
 As we saw in the previous section, there are many integral 
 representations of the solutions of the wave equation on $[0,1]^{2}$.
 
 Most of them  appear as 
 $$F(M)= \int F(S)\ V_{0}(M,dS)+ \int F(S)\ V_{1}(M,dS),$$
  where 
  $V_{0}(M,dS)$ is a continuous family of bounded  measures whose 
  support is the interval  $[M_{-},M_{+}]$ on the boundary $\{y=0\}$, and $V_{1}(M,dS)$ is a 
  continuous family of bounded measures with support $\Delta_{M}$.
  
  In general,  those representations lead  to a 
  unique representation
  $$F(M)= \int_{[M_{-},M_{+}]} F(S)\ W(M,dS),$$ for a continuous 
  family of bounded measures with support in $[M_{-},M_{+}]$. The 
  crucial point is that in some situations the measure $W(M,dS)$ may be 
  positive even if $V_{1}(M,dS)$ is not. In the case of Achour--Trim\`eche's 
  theorem however, the measure $W$ is positive only on some symmetric 
  functions.
   
  To understand these representations, we shall consider a more 
  general setting.
 
  Consider a separable compact Hausdorff space $E$, and two continuous 
  families $V_{0}(M,ds)$ and $V_{1}(M,ds)$ on $E$ (two kernels). We shall identify 
  such a family with the operator 
  $$F\mapsto V_{i}(F)(M)= \int_{E} F(y)\ V_{i}(M,dS),$$ which maps the 
  Banach space $C(E)$ of continuous functions into itself. The 
  identity operator corresponds to the kernel $I(M,dS)= 
  \delta_{M}(dS)$, and the 
  composition of kernels
  $$V\odot W (M,dS) = \int _{U} \ V(M,dU)W(U,dS)$$ corresponds then 
  to the operator composition.
  
  We set $$\left\|V\right\|= \sup_{M}\int_{E}\left|V(M,dS)\right|,$$
  which is the operator norm.
  
  Then we have
  \begin{lemma}
  Consider  some continuous function $F\in C(E)$ satisfies
  \begin{equation}\label{eqgal}\forall M\in E,~ F(M) = \int F(S)\ V_{0}(M,dS) + \int_{E} F(S) \
  V_{1}(M,dS).\end{equation}
  
  If the series 
  \begin{equation}\label{condserie}\sum_{n} \left\|V^{\odot n}_{1}\right\|
   \end{equation} converges,
  then setting 
  $${\cal E}(V_{1})=W_{1}= \sum_{n} V_{1}^{\odot n},$$
  one has 
  $$F= W_{1}\odot V_{0}(F).$$
  
  In particular, if $V_{0}$ is supported by some closed subset 
  $E_{0}$ of $E$, and if the condition \eqref{condserie} is satisfied 
  for $V_{1}$, then there exists a unique solution $F$ to the 
  equation 
  \eqref{eqgal} given the restriction $F_{0}$ of $F$ on $E_{0}$. 
  
  Moreover, if the kernels $V_{0}$ and $V_{1}$ are non-negative and if 
  $F_{0}$ is non-negative, so is $F$.
  \end{lemma}
  
  \begin{proof}
  The proof is straightforward and is just the classical 
  representation of $(I-V_{1})^{-1}$ as $\sum_{n} V_{1}^{\odot n}$. 
  Observe moreover that if $V_{1}$ is non-negative, so is $W_{1}= 
  {\cal E}(V_{1})$, and that the representation of the solution is 
  then 
  given by a non-negative kernel.\qed
\end{proof}

 In what follows, and to apply these lemmas,  we shall consider the 
 case where $E$ is a compact subset of $\R\times [0,\infty[$, and 
 where, for some point $M\in E$, the measure $V_{0}(M,dS)$ is supported by 
 $[M_{-},M_{+}]$,  $V_{1}(M,dS)$ is supported by 
 $\Delta_{M}$, and has a bounded density $a(M,S)$ with respect to 
 the Lebesgue measure on the product, in which case we write
 
 \begin{equation}
 \label{eqint2}
 F(M)= \int_{[M_{-},M_{+}]} F(S)\ V_{0}(M,dS)+ \int_{\Delta_{M}} F(S) 
 a(M,S)\ dS.
  \end{equation}
 
 Then we have 
 \begin{proposition} \label{repr3}
 Let $\kappa$ a uniform bound on 
 $\left|a(S,U)\right|$, $(S,U) \in \Delta_{M}$ in \eqref{eqint2}.
 Then, on $\Delta_{M}$,  we have 
 $$\left\|V_{1}^{\odot n}\right\|\leq 
 \frac{\kappa^{n}\left|\Delta_{M}\right|^{n}}{(n!)^{2}},$$
 where $\left|\Delta_{M}\right|$ denotes the area of the triangle 
 $\Delta_{M}$.
 \end{proposition}
 
 \begin{proof}
 Recall the partial order $(S_{1} \leq S_{2}) \iff S_{1} 
 \in \Delta_{S_{2}}$.

 Then one has 
 $$V_{1}^{\odot n}(M,dS)= \Ind_{\Delta_{M}}(S)a_{n}(M,S) dS,$$
  where 
  \begin{multline*}
a_{n}(M,S)=\\ \int_{S\leq S_{1}\leq \cdots\leq S_{n-1}\leq M} 
  a(M,S_{n-1})a(S_{n-1},S_{n-2})\ldots a(S_{1},S)\ dS_{1}\ldots\  
  dS_{n-1}.
\end{multline*}
  
  It is easy to see by induction that 
  $$ \left|a_{n}(M,S)\right| \leq \kappa^{n}\frac{ 
  \left|[S,M]\right|^{n}}{n!^{2}},$$
  where
  $\left|[S,M]\right|$ denotes the area of the rectangle $\{U~\mid ~ 
 S\leq U\leq M\}$.
 
 The conclusion follows easily from this estimate. \qed
\end{proof}

 Considering the representation  of the 
 solutions of the wave equation given in  Proposition \ref{repr2}, we 
 see that any  continuous solution may be represented as
 $$F(M)= \int_{[M_{-}M_{+}]} F(S)\ V(M,dS)$$ where 
 $V(M,dS)$ has two Dirac masses at the points $M_{-}$ and $M_{+}$ and 
 has a bounded density on $(M_{-},M_{+})$. The smoothness of this 
 density depends of course on the smoothness of the function $a$ 
 itself (and may be analyzed through the convergence of the series 
 that we just described). For example, if $a$ has $k$ bounded 
 derivative, then so has the density.

  We may observe the following.
  \begin{corollary}\label{comp1} Consider a solution $F$ of Equation \eqref{eqgal}.
  Assume that $V_{0}$ is non-negative 
  and supported by 
  $E_{0}\subset E$, and that $V_{1}\geq V_{2}$, in the sense that 
  $V_{3}= V_{1}-V_{2}$ is a non-negative kernel. Suppose that $V_2$ 
  and $V_{3}$ satisfy the growth condition \eqref{condserie} and moreover that 
  ${\cal E}(V_{2})= W_{2}$ is non-negative.
  
  Then, if the restriction  $F_{0}$ of $F$ on $E_{0}$ is non-negative, then $F$ is non-negative everywhere on $E$.
 
 \end{corollary}
 
 \begin{proof} Once again, this is straightforward. Setting $W_{2}= 
 (I-V_{2})^{-1}$, one has
 $$F= W_{2}\odot V_{0} (F)+ W_{2}\odot V_{3}(F),$$
 which is an equation of the same type, but with non-negative kernels.\qed
\end{proof}

We may apply this for example for the solutions of \eqref{eqint2}~:

\begin{corollary} If $V_{0}(M,dS)$ is non-negative and $a(U,S)\geq -C$, on $S\leq 
U\leq M$, where $$C= \frac{\mu_{0}^{2}}{2\left|\Delta_{M}\right|}$$ 
and $\mu_{0}$ is the first $0$ of the Bessel $G$ function which is 
solution on $(0, \infty)$ of 
$$ G''+ \frac{G'}{x} = -G,~ G(0)=1, ~ G'(0)=0,$$
then any continuous solution of Equation \eqref{eqint2} which is non-negative on the boundary $\{y=0\}$ is non-negative on $\Delta_{M}$.
\end{corollary}

\begin{proof}
It is a simple application of the previous Corollary \ref{comp1} with 
$V_{2}(M,dS)=-C\Ind_{\Delta_{M}}(S)dS$.

In this case, it is not hard to see that 
$${\cal E}(V_{2})(M,dS)= \Ind_{\Delta_{M}}(S) F(\left|[S,M]\right|)\ dS,$$
where $\left|[S,M]\right|$ denotes the Lebesgue measure of the 
rectangle $[S,M]=\{U~\mid ~ S\leq U\leq M\}$ and 
$$F(x) = \sum_{n} \frac{(-Cx)^{n}}{(n!)^{2} }.$$
The function $F$ is the solution of 
$$xF''+ F'= -CF, ~ F(0)=1, F'(0)=-C,$$
which is related to the function $G$ through the change of variable 
$x={z^{2}}/{4C}$.
The function $G$ is non-negative on $[0, \mu_{0})$ and this gives the 
result provided one observes that $\left|[S,M]\right|\leq 
\frac{\left|\Delta_{M}\right|}{2}$.\qed
\end{proof}

\begin{rk}
One may also observe that $
{1}/{\mu_{0}^{2}}$ is the fundamental 
eigenvalue of the Laplace operator on the unit ball of $\R^{2}$ with 
the Dirichlet boundary conditions, the function 
$G({\left\|x\right\|}/{\mu_{0}})$ being the corresponding 
eigenvector. 
\end{rk}

All these considerations  provide many criteria on the 
function $a$ such that the associated Neumann basis on $I$ has the 
hypergroup property at the left end point of the interval. In the 
next section, we shall deal with Gasper's theorem, where $I= [0,\pi/2]$ 
and $a(x)= \alpha\tan x - \beta\cot x$, with $\alpha\geq \beta>-1$. The reader should 
check that no one of these criteria may apply on this example. 
Achour--Trim\`eche's theorem shows that the hypergroup property holds for this 
example in the symmetric case $\alpha=\beta$, even on any symmetric  
(around $\pi/4$) 
subinterval of $[0,\pi/2]$. But we do not even know for the moment if 
the hypergroup property holds in the general case on any symmetric 
subinterval of $[0,\pi/2]$. (Yet it is true for small subintervals 
and also provided that the parameters $\alpha$ and $\beta$ belong to 
some specific domains that we shall not describe here).

\section{The case of Jacobi polynomials: Gasper's 
theorem}\label{Gasper}

Gasper's theorem states the hypergroup property for the family of 
Jacobi polynomials.
The case of Jacobi polynomials may be considered as a special case of 
a
Sturm--Liouville basis on $[0, \pi/2]$. In this situation, both the GKS and the HGP 
property hold \cite{Gasper70,Gasper71,Gasper72}. Actually, it is a unique situation for 
orthogonal polynomials, since they are the only ones, up to a linear change of variables, 
for which the HGP property holds (see \cite{CS90,CMS91,CS95}) (under 
some mild extra condition on the support of the measure which 
represents the product formula). In the case of 
symmetric Jacobi polynomials (known as Gegenbauer or ultraspherical 
polynomials), the HGP property may be seen as a particular example of 
Achour--Trim\`eche's theorem (although in this case the measure has a density 
which vanishes on the boundary). But in the general case, as we 
already mentioned, none of the extension we gave of Achour--Trim\`eche's theorem covers 
this result. Even worse, we do not know if the HGP property holds for 
any symmetric subinterval of $[0,\pi/2]$. 

The Jacobi polynomials are a quite universal object, since they are 
basically the unique examples of a family of orthogonal polynomials 
which are also eigenvectors of Sturm--Liouville operators (together 
with their limiting cases the Hermite and Laguerre polynomials, see 
\cite{Mazet98-1}). On the other hand, for special values of the 
parameters, they may be considered as eigenvectors of rank-one 
symmetric compact spaces (here, with our notations, it is for the 
parameters $(1,p)$, $(2,p)$, $(4,p)$ and $(p,p)$, with $p\in \N$). But for a wider range of parameters 
($p,q$ $\in \N$ ), they may 
be seen as eigenvectors of a Laplace operator on a $p+q-1$ 
dimensional sphere. The special case where $p=q$ is much simpler, 
since then one may consider a $p$-dimensional sphere.

In this section, after a short introduction on Jacobi polynomials and the statement of  the hypergroup property for these polynomials, we 
present the simpler case of symmetric Jacobi polynomials, where the 
convolution structure has a nice geometric interpretation for $p\in 
\N$. This interpretation is for example described in \cite{Bingham}. For the dissymmetric case, although the Jacobi polynomials still  have a simple geometric 
interpretation too when the 
parameters are integers, the convolution structure is far less obvious.

\subsection{Jacobi polynomials}\label{JacobiPolyn}

This polynomial family  is defined for some positive parameters $p$ and $q$ as  the family of 
orthogonal polynomials associated with
the measure 
$$\mu_{p,q}(dx) = C_{p,q}(1-x)^{\frac{q-2}{2}}(1+x)^{\frac{p-2}{2}} dx $$ 
on $[-1,1]$, $C_{p,q}$ being a normalizing constant such that 
$\mu_{p,q}$ is a probability measure.

These polynomials are also the eigenvectors of the operator 
$$L_{p,q}f(x) = (1-x^{2}) f''(x) - 
\left(q\frac{x+1}{2}+p\frac{x-1}{2}\right)f'(x)$$ on $[-1,1]$. 
If $P_{k}^{p,q}$ is the polynomial of degree $k$, one has 
$$L_{p,q}P_{k}^{p,q}= -k\left(\frac{p+q}{2} +k-1\right) P_{k}^{p,q}.$$

\begin{rk}
These polynomials are traditionally parametrized by $\alpha=\frac{q-2}{2}$ and  
$\beta=\frac{p-2}{2}$ with $\alpha,\beta > -1$, from \cite{Szego75} or \cite{Gasper70,Gasper71,Gasper72}.
\end{rk}

If we change $x= \cos(2\theta)$, $\theta\in [0,\frac{\pi}{2}]$, then 
this operator is turned into
$$L_{p,q}f(\theta)= \frac{1}{4}\left[ f''(\theta) + \left((q-1)\cot (\theta) 
-(p-1)\tan (\theta)\right)f'(\theta)\right].$$ 

We see then that Jacobi polynomials is one example of a Neumann basis associated with a Sturm--Liouville operator
(except that the density of the measure vanishes on the boundary 
points, for parameters larger than 2). We may also observe that the measure is log-concave as soon 
as $p$ and $q$ are in $[1,\infty)$.

When $p$ and $q$ are integers, one may see the operator $L_{p,q}$ as 
the action of some spherical laplacian on a quotient of the sphere.

More explicitly,  we set $N= p+q$. Let us denote by 
$\left|X\right|$ the euclidean norm  of a point $X$ in $\R^{N}$, and let 
$\sphere^{N-1}$ be the unit sphere. We consider the Laplace 
operator $\Delta_{\sphere}$  on the unit sphere $\sphere^{N-1}$ in $\R^{N}$: this is the 
restriction to the sphere of the usual Laplace operator on $\R^{N}$ 
acting on function which are defined in a neighborhood of the sphere 
and do not depend on the radius of the point.

We parametrize $\sphere^{N-1}$ as 
\begin{equation}\label{coordspheres}X= \left(\sqrt{\frac{1+x}{2}}X_{1}, \sqrt{\frac{1-x}{2}}X_{2}\right),
\end{equation}
where 
$X_{1}\in \sphere^{p-1}$ 
$X_{2}\in \sphere^{q-1}$, and $x\in [-1,1]$. The
action of $\Delta_{\sphere}$ on a function which depends only on 
$x$ gives again a function of $x$, and 
we have
$$\Delta_{\sphere}(h)(x)=  4L_{p,q} 
(h)(x).$$

We shall say that such a function on the sphere which depends only on 
$x$ (that is which depends only on the norm of the projection of $X$ 
onto $\R^{p}$) has the invariance 
$SO(p)\times SO(q)$, where the action of $SO(p)\times SO(q)$ is 
obtained by the action of the first component on $X_{1}$ and of the 
second on $X_{2}$.

The measure  $\mu_{p,q}$ is the invariant measure for the operator $L_{p,q}$ 
and the uniform measure is the invariant measure for the Laplace 
operator on the sphere. This shows that $\mu_{p,q}$ is the image of 
the uniform measure on the sphere (normalized as to be a probability 
measure) under the map $X\mapsto x$ of Formula \eqref{coordspheres}.

In fact, under this map $X\mapsto (X_{1},X_{2},x)$, it is straightforward to see that the 
uniform measure $\sigma_{p+q-1}$ on $\sphere^{p+q-1}$ is transformed into 
$\sigma_{p-1}\otimes\sigma_{q-1}\otimes\mu_{p,q}$.

Thanks to this remark, consider $N\geq p$ and look at the projection 
$\pi(X)$ from $\sphere^{N-1}$ onto the unit ball in $\R^{p}$. (That is the 
orthogonal projection when the sphere is imbedded into $\R^{N}$). If we 
set $x= 2\left|\pi(X)\right|-1$ and $X_{1}=\frac{ 
\pi(X)}{\left|\pi(X)\right|}$, we  see that the image measure of 
$\sigma_{N-1}$ under $X\mapsto (x,X_{1})$ is $\mu_{p,N-p}\otimes 
\sigma_{p-1}$. This remark shall be used in Paragraph 
\ref{JacobDissym}.

\subsection{Gasper's result}
Gasper proved the following product formula which gives the HGP property for Jacobi polynomials, applying Proposition \ref{HGP2}.
\begin{theorem}[Gasper]
Let $p,q>0$ and $-1< x, y <1$. Then 
\begin{itemize}
\item
we have the following product formula:
$$
\forall k, \quad
\frac{P_k^{p,q}(x)P_k^{p,q}(y)} {P_k^{p,q}(1)}=
\int P_k^{p,q}(z)\ m_{p,q}(x,y,d{z})$$
where $m_{p,q}(x,y,dz)$ is a Borel measure on $[-1;1]$;
\item
the measure $m_{p,q}$ is positive (and then is a probability measure) if and only if 
$$(p,q)\in \{q \geq p\} \cap \{p \geq 1 \text{ or } p+q \geq 4\};$$
\item
moreover, if $q>p>1$, $m_{p,q}(x,y,dz)$ is absolutely continuous with respect to $\mu_{p,q}$, with density in $L^2(\mu_{p,q})$, so that
$$K_{p,q}(x,y,z)=\sum_k \frac{P_k^{p,q}(x)P_k^{p,q}(y)P_k^{p,q}(z)}{P_k^{p,q}(1)} \geq 0,$$
with convergence of the sum for almost every $z$.
\end{itemize}
\end{theorem}
The original Gasper's proof (see \cite{Gasper71,Gasper72}) consisted in the explicit computation of the sum $K_{p,q}(x,y,z)$ using formulae on special functions like Bessel's and hypergeometric functions.
There had been many other proofs of this property. For example,  Koornwinder derived it in
\cite{Koorn73} from the addition formula of Jacobi polynomials and he found an other proof in \cite{Koorn74}, that we discuss next.

Here we restrict ourself to prove the HGP property in the symmetric case ($p=q$) 
and in the case when $q> p>1$. In the latter case, we follow a proof given by Koornwinder
in \cite{Koorn74}. However, his argument was based on an integral representation formula 
of the polynomials (our Lemma \ref{lem2}), whose proof, as we found in literature
(see \cite{Askey74} together with \cite{AskeyFitch69}), relies on computational considerations
on hypergeometric functions. In Section \ref{JacobDissym}, we shall give a more geometric 
interpretation of this formula, at least when $p$ and $q$ are integers (it appears finally that
the interpretation of Jacobi polynomials as harmonic functions was already known 
-- see \cite{braak-meul,Koorn71,Koorn73} -- but it seems that it was not yet directly used to 
derive Koornwinder's representation formula).

\subsection{The special case of ultraspherical 
polynomials ($p=q\geq1$)}

 In the case when $p=q$, Jacobi polynomials are called ultraspherical 
 polynomials. In this case, there is a much simpler representation of 
 $L_{p,p}$
 when $p\in \N^{*}$ as the action of the Laplace operator on the sphere 
 $\sphere^{p}$ (and not on $\sphere^{2p-1}$ as before).
 
 When $p$ is a positive integer, then the hypergroup property has a 
 simple geometric interpretation, and thus the property is quite 
 easy  to establish. This easily extends to the case when $p\notin 
 \N$, by a simple extension of the formulae. This is what we are 
 going to see in this paragraph.
 
 Consider a smooth function $F~:~\sphere^{p}\mapsto \R$ which depends 
 only on the first coordinate. To fix the ideas, let $F(X)= f(X\cdot 
 e_{1})$, where $e_{1}$ is the first unit vector in $\R^{p+1}$, and 
 $Y\cdot X$ denotes the standard scalar product in $\R^{p+1}$. Then, 
 if $\Delta_{\sphere^{p}}$ is the Laplace operator on $\sphere^{p}$, we have
 $$\Delta_{\sphere^{p}} = L_{p,p}(f)(X\cdot e_{1}).$$
 
 As before, the image measure of the uniform measure on the sphere 
 through the map $X\mapsto x=X\cdot e_{1}$ is the invariant measure 
 for $L_{p,p}$, that is $\mu_{p,p}$. Moreover, we may parametrize 
 $\sphere^{p*}=\sphere^{p}\setminus \{e_{1},-e_{1}\}$ by
 \begin{equation} \label{paramsphere} X= \left(x, \sqrt{1-x^{2}}X_{1}\right),\end{equation} where $x\in (-1,1)$ is the first 
 coordinate of the point $x\in \sphere^{p}\subset \R^{p+1}$, and $X_{1}\in 
 \sphere^{p-1}$. Through this map $\sphere^{p*}\mapsto (-1,1)\times \sphere^{p-1}$, the image measure of $\sigma_{p}$ is 
 $\mu_{p,p}\otimes \sigma_{p-1}$.

From that, we see that if $P_{k}^{p,p}$ is the $k$-th ultraspherical polynomial, 
 then $P_{k}^{p,p}(X.e_{1})$ is an eigenvector of $\Delta_{\sphere^{p}}$, 
 with eigenvalue $\lambda_{k}^{p}= -k(k+p-1)$. 
 
Observe that for any point $Y$ on the sphere, $P_{k}^{p,p}(Y.X)$ is 
again an eigenvector on the sphere with the same eigenvalue. (This 
comes from the fact that the Laplace operator on the sphere is 
invariant under rotations.)

Now, if we take two points $Y$ and $Z$ on the sphere, 
$P_{k}^{p,p}(Y.X)$ and $P_{k}^{p,p}(Z.X)$ are two eigenvectors of 
$\Delta_{\sphere^{p}}$, with the same eigenvalue. Let us compute their 
scalar product in $L^{2}(\sphere^{p})$

$$H(Y,Z)= \int_{\sphere^{p}} P_{k}^{p,p}(Y.X)P_{k}^{p,p}(Z.X)\ 
\sigma_{p}(dX).$$

Obviously, $H(Y,Z)$ is a smooth function, taking values in $[-1,1]$, 
and if $R$ is any rotation, $H(Y,Z)= H(RY,RZ)$. From this we see that
$H(Y,Z)= h(Y\cdot Z)$, and we may write the function $h$ in terms of 
ultraspherical polynomials
$$h= \sum_{r} a_{r} P_{r}^{p,p}.$$

We have
\begin{eqnarray*} a_{r}&=& \int h(x)P_{r}(x)\ \mu_{p,p}(dx)= \int_{\sphere^{p}} 
h(Y.Z)P_{r}^{p,p}(Y.Z)\ \sigma_{p}(dZ)\\&=& \int_{\sphere^{p}}\int_{\sphere^{p}} 
P_{r}^{p,p}(Y.Z)P_{k}^{p,p}(Y.X)P_{k}^{p,p}(Z.X) \ 
\sigma_{p}(dX)\sigma_{p}(dZ).\end{eqnarray*}
Using Fubini's theorem and the orthogonality of eigenvectors 
associated with different eigenvalues, we see that $a_{r}=0$ unless 
$r=k$. We therefore see that
$$\int_{\sphere^{p}} P_{k}^{p,p}(Y.X)P_{k}^{p,p}(Z.X)\ \sigma_{p}(dX)= 
a_{k}P_{k}^{p,p}(Y.Z).$$
To compute $a_{k}$ we choose $Y=Z$, from which we get, if we remember 
that the polynomials $P_{k}^{p,p}$ have norm $1$ in $L^{^2}$, that
$$a_{k}P_{k}^{p,p}(1)=1.$$

Now, if we rewrite this formula for $Y= e_{1}$ and through the 
parametrization described above in \eqref{paramsphere}, then we get, for
$Z=(z, \sqrt{1-z^{2}}Z_{1})$,
$$\frac{P_{k}^{p,p}(z)}{P_{k}^{p,p}(1)} = \int_{x,X_{1}} 
P_{k}^{p,p}(x)P_{k}^{p,p}(zx+ \sqrt{1-z^{2}}\sqrt{1-x^{2}} 
X_{1}Z_{1})\ \mu_{p,p}(dx)\sigma_{p-1}(dX),$$
while, for $k\neq l$,
$$\int_{x,X_{1}} 
P_{k}^{p,p}(x)P_{l}^{p,p}(zx+ \sqrt{1-z^{2}}\sqrt{1-x^{2}} 
X_{1}Z_{1})\ \mu_{p,p}(dx)\sigma_{p-1}(dX)=0.$$
This may be rewritten as
$$\int_{x,t} 
P_{k}^{p,p}(x)P_{l}^{p,p}(zx+ \sqrt{1-z^{2}}\sqrt{1-x^{2}} 
t)\ \mu_{p,p}(dx)\mu_{p-1,p-1}(dt)=\delta_{kl}\frac{P_{k}^{p,p}(z)}{P_{k}^{p,p}(1)}.$$

This last formula may be turned into an explicit representation
$$\int_{x,t} 
P_{k}^{p,p}(x)P_{l}^{p,p}(y)\ k_{p}(z,dx,dy)=\delta_{kl}\frac{P_{k}^{p,p}(z)}{P_{k}^{p,p}(1)}$$
 for some probability 
kernel $k_{p}(z,dx,dy)$, which gives the hypergroup property thanks 
to Proposition \ref{HGP3}.

When $p$ is not an integer, since we have an explicit representation 
of the kernel $k_{p}(z,dx,dy)= K_{p}(x,y,z) \mu(dx)\mu(dy)$, it is a 
simple verification to check that the function
$K_{p}(x,y,z)$ satisfies $L_{x} K_{p}= L_{y} K_{p}$ together with
$$\lim_{y\to 1}K_{p}(x,y,z)\ \mu_{p,p}(dz)=\delta_{x}(dz),$$ which is 
enough to get the HGP property at the point $1$.

Moreover, the convolution associated with this hypergroup structure 
is quite easy to understand when $p$ is an integer. 

Let us say that a probability measure $\mu$ on the sphere $\sphere^{p}$ is 
zonal around $X\in \sphere^{p}$ if it is invariant under any rotation $R\in 
SO(p+1)$ such that $RX=X$.

Given any probability measure $\mu$ on $[-1,1]$, and any $X\in 
\sphere^{p}$, we may lift $\mu$ into a unique probability measure $\hat \mu$ 
which is zonal around $X$ such that the image measure of $\hat \mu$ 
under the projection $\pi(Y)= Y\cdot X$ from $ \sphere^{p}$ onto $ [-1,1]$  is 
$\mu$.

Now, let us choose $e_{1}\in \sphere^{p}$, and consider  two probability measures $\nu_{1}$ and $\nu_{2}$ on 
$[-1,1]$. We may lift $\nu_{1}$ into a probability measure $\hat 
\nu_{1}$ on $\sphere^{p}$, which is zonal around $e_{1}$. Then, we choose a  
random
point in $\sphere^{p}$ according to $\hat \nu_{1}$. Then, given $X$, we 
consider the lift of $\nu_{2}$ which is zonal around $X$ and choose a 
random point $Y$ according to this measure. Then, the resulting 
law of $Y$ is zonal around $e_{1}$, and we project this measure into 
a new measure $\nu_{1}*\nu_{2}$. It is an exercise to show that this 
convolution is the convolution associated with the hypergroup structure 
in this case.

\subsection{The case of dissymmetric Jacobi polynomials ($q>p>1$)}\label{JacobDissym}

Although the dissymmetric Jacobi polynomials may be interpreted as 
eigenvectors of the Laplace operator on the sphere $\sphere^{p+q-1}$, it is 
far from trivial to prove the hypergroup property even in the case 
where $p$ and $q$ are integers. Nevertheless, the proof that we 
present below for  completeness and which is due 
essentially to Koornwinder \cite{Koorn74} has also some simple 
interpretation when $p$ and $q$ are integers in terms of harmonic 
analysis in $\R^{p+q}$.

Koornwinder's proof relies   on two facts, given in the following 
Lemmas \ref{lem1} and \ref{lem2}. In what follows, and to lighten the notations, we 
remove the indices $p$ and $q$ from the definitions of the 
polynomials $P_{k}^{p,q}$.

\begin{lemma} (Bateman's formula)\label{lem1} Let $b_{k,r}$ the the coefficients such that
$$\frac{P_{k}(s)}{P_{k}(1)}= \sum_{r=0}^{k} b_{k,r} (s+1)^{r}.$$
Then
$$\frac{P_{k}(s)P_{k}(t)}{P_{k}(1)^{2}}= 
\sum_{r=0}^{k}b_{k,r}\frac{(s+t)^{r}}{P_{r}(1)} 
P_{r}\left(\frac{1+st}{s+t}\right).$$

\end{lemma}

\begin{lemma}\label{lem2}(Koornwinder's formula)

\begin{multline*}
\frac{P_{k}(x)}{P_{k}(1)} = \int_{[-1,1]^{2}}\Bigg[ \frac{
2(1+x)-(1-x)(1+u)}{4}\\
+\frac{\I\sqrt{2}\sqrt{1-x^{2}}\sqrt{1+u}\,v}{4}\Bigg]^{k}
\mu_{p,q-p}(du)\mu_{p-1,p-1}(dv).$$
\end{multline*}
\end{lemma}

Before going further, let us show that this implies the HGP property 
at the point $x_{0}=1$.  In fact, we shall use the characterization of 
the hypergroup property given by Proposition \ref{HGP2}.

For that, we replace in Bateman's formula of \ref{lem1} the 
representation given by Koornwinder's formula. For this, we observe 
that, if $(s,t)\in  (0,1)^{2}$, then 
$$\left(\frac{1+st}{s+t}\right)^{2}>1,$$
and therefore if we set $x={(1+st)}/{(s+t)}$, we may replace  by 
analytic continuation
$\I{}\sqrt{1-x^{2}}$ by $\sqrt{x^{2}-1}$.

Then, if we set
$$\psi(s,t,u,v)= \frac{2(1+x)-(1-x)(1-u)+ 
\sqrt{2}\sqrt{x^{2}-1}\sqrt{1+u}\,v}{4},$$
one has
$$\frac{P_{k}(s)P_{k}(t)}{P_{k}(1)^{2}}= \int \sum_{r=0}^{k} b_{k,r} 
[(s+t)\psi(s,t,u,v)]^{k}\ \mu_{p,q-p}(du)\mu_{p-1,p-1}(dv).$$

From the definition of the coefficients $b_{k,r}$, we get then
$$\frac{P_{k}(s)P_{k}(t)}{P_{k}(1)^{2}}= \int 
\frac{P_{k}\Big((s+t)\psi(s,t,u,v)-1\Big)}{P_{k}(1)}\ 
\mu_{p,q-p}(du)\mu_{p-1,p-1}(dv).$$
If we define $m_{p,q}(s,t,dz)$ to be the image measure of 
$\mu_{p,q-p}(du)\mu_{p-1,p-1}(dv)$ under the map
$$(u,v)\mapsto (s+t)\psi(s,t,u,v)-1,$$ one gets
$$\frac{P_{k}(s)P_{k}(t)}{P_{k}(1)}= \int P_{k}(z)\ m_{p,q}(s,t,dz),$$ which 
is the announced result.

Of course, one has to check that the image measure is indeed 
supported by $[-1,1]$, but this point is left to the reader.

We now give the proof of  Lemmas \ref{lem1} and \ref{lem2}. As it shall turn out, they 
rely on  elementary considerations on the interpretations of the 
operator $L_{p,q}$. For the moment, we restrict ourselves to the case 
where $p$ and $q$ are positive integers, and we shall interpret those 
formulae in term of the Laplace operator on $\R^{p+q}$.

First observe that, given any function $f$ on $[-1,1]$, we may lift 
this function on the sphere $\sphere^{p+q-1}$ into a function which has the 
$SO(p)\times SO(q)$ invariance. Namely, using the parametrization of 
the sphere given in \eqref{coordspheres}, we set
$$F\left(\sqrt{\frac{1+x}{2}}X_{1}, \sqrt{\frac{1-x}{2}}X_{2}\right)= f(x).$$
Notice that in this formula,  
$$x=\left|\pi_{1}X\right|^{2}-\left|\pi_{2}(X)\right|^{2},$$ where 
$\pi_{1}$ and $\pi_{2}$ are the orthogonal projections on $\R^{p}$ 
and $\R^{q}$ when the sphere is imbedded into $\R^{p+q}$.  Let us call 
$U(f)$ such a lift of  a function from $[-1,1]$ onto the sphere.

Now, if $P_{k}$ is the Jacobi polynomial of degree $k$, the 
corresponding function $F_{k} = U(P_{k})$ is an eigenvector of the Laplace 
operator on the sphere, and therefore the restriction to the sphere of a 
harmonic polynomial of degree $2k$. Therefore, if we parametrize a 
point $Z$ in $\R^{p+q}$ by $R= \left|Z\right|^{2}$ and $X= 
\frac{Z}{\left|Z\right|} $, we may see that the function 
$R^{k}F_{k}(X)$ is harmonic in $  \R^{p+q}$.

This may be seen in another way as we may write the Laplace 
operator in those coordinates 
$$\Delta= \frac{\partial^{2}}{\partial_{R}^{2}}+ 
\frac{N}{2R}\frac{\partial}{\partial_{R}} + 
\frac{1}{4R^{2}}\Delta_{\sphere},$$
where $N=p+q$ and  $\Delta_{\sphere}$ is the Laplace operator on $\sphere^{N-1}$. 
(It does not look as usual because of the change of $r= \left|x\right|$ 
into $r^{2}=R$.)
 Since 
 $$\Delta_{\sphere} F_{k}= 4U\left(L_{p,q}P_{k}\right)= -4k\left(k+\frac{N}{2}-1\right)F_{k},$$
  one may check directly that $H(R,X)= R^{k}F_{k}(X)$ is a solution of 
  $\Delta H=0$.
  
  In other words, the solutions of 
  $$\left(\frac{\partial^{2}}{\partial_{R}^{2}}+ 
\frac{N}{2R}\frac{\partial}{\partial_{R}} + 
\frac{1}{R^{2}}L_{p,q}\right)F=0$$ correspond to harmonic functions in 
$\R^{p}\times \R^{q}$ which are radial in both components (bi-radial 
harmonic functions). If $(X,Y)$ are the two component of a point in 
$\R^{p+q}$, then this harmonic function is
$$(\left|X\right|^2+\left|Y\right|^{2})^{k}
P_{k}\left(\frac{\left|X\right|^2-\left|Y\right|^{2}}{\left|X\right|^2+\left|Y\right|^{2}}\right).$$

\begin{proof} (Of Bateman's formula \ref{lem1}.)

Let $L= L_{p,q}$. The function $K(s,t)=P_{k}(s)P_{k}(t)$ is a solution of 
the wave equation 
$$(L_{s}-L_{t})K=0.$$
In order to prove the assertion, which amounts to verify the identity 
of two polynomials, it is enough to check it on an open set. We shall 
choose to prove it on the set $\{s\in (-1,1), ~ t>1\}$, on which the  wave 
equation $(L_{s}-L_{t})K=0$ becomes an elliptic equation.

On the other hand, consider a solution $G(R,x)$ of 
$$\left(\frac{\partial^{2}}{\partial_{R}^{2}}+ 
\frac{N}{2R}\frac{\partial}{\partial_{R}} + 
\frac{1}{R^{2}}L_{p,q}\right)G=0,$$
and perform the change of variable
$$R=s+t; ~ x= \frac{1+st}{s+t}.$$
This equation becomes $(L_{s}-L_{t})G=0$. (We shall leave the 
computation to the reader, since it is just brute  calculus.)

This strange (and miraculous) change of variables may be much 
understood if we first operate a change of variables to reduce the 
leading terms to $\partial_{x}^{2}+ \partial_{y}^{2}$ in both 
equations, and then observe that the transformation we made is 
conformal in $\R^{2}$, and thus preserves the leading terms. But we 
could find no simple geometric transformation, even in the case where $p$ 
and $q$ are positive integers, to understand this change of a 
bi-Jacobi equation into a bi-radial harmonic function.

Therefore, the right-hand side of Bateman's formula is a solution of 
the wave equation $(L_{s}-L_{t})F=0$. The coefficients $b_{k,r}$ are 
computed in such a way that the two polynomials coincide on $t=1$.

To see that they must coincide everywhere, it is enough to remark 
that if two polynomials $A(s,t)$ and $B(s,t)$ in $(s,t)$ are 
solutions of the wave equation which coincide on $t=1$, they coincide 
everywhere. Indeed, we may write
$$A(s,t)= \sum_{r=0}^{k} a_{r}P_{r}(s)P_{r}(t), ~B(s,t)= 
\sum_{r=0}^{k}b_{r}P_{k}(s)P_{k}(t),$$
and identifying the values in $t=1$ produces $a_r= b_{r}$, 
$r=0,\ldots,k$.\qed
\end{proof}

We now turn to the proof of Koornwinder's formula.

We begin with a lemma.
Here, we shall use for the first time that $q>p$.  

\begin{lemma} \label{repreOq}
Consider a 
bi-radial analytic function $H$ on $\R^{p+q}$,  that is to say
$H(X,Y)= h\big(\left|X\right|^{2}, \left|Y\right|^{2}\big)$ where $h$ is an 
analytic function on $\R^{2}$. If moreover $H$ is harmonic on $\R^{p+q}$, 
then it holds
\begin{equation} \label{reprebirad} 
H(X,Y)= \int_{S0(q)} H(X+\I{}\pi(RY),0)\ \nu(dR)
\end{equation}
 where $\nu(dR)$ is the Haar measure on the group $SO(q)$, and $\pi$ 
 is  the orthogonal projection  from $\R^{q}$ onto $\R^{p}$.

 As a consequence, if $H(X,Y)= h(\left|X\right|^{2}, 
 \left|Y\right|^{2})$ and $g(x)= h(x,0)$, then
\begin{multline}\label{reprOq2}
{h(\left|X\right|^{2},\left|Y\right|^{2})=} \\
 \int_{[-1,1]^{2}} g
\bigg(\left|X\right|^{2}-\frac{1+u}{2} \left|Y\right|^{2}+\I{}\sqrt{2}\left|X\right|\left|Y\right|\sqrt{1+u}s\bigg)
\  d\mu_{p,q-p}(u)d\mu_{p-1,p-1}(s).
\end{multline}

 \end{lemma}
 
 In practice, we shall just apply this lemma with polynomials 
 functions $h$.
 
 \begin{proof} (Of Lemma \ref{repreOq}.)
 
 The proof comes from the following remark. We observe that if 
 $F$ is any analytic radial function in $\R^{p}$, namely $F(X)=f(\left|X\right|^{2})$ where $f$ is real analytic, then 
 $F(X+\I{}Y)$ is a solution in $\R^{p}\times \R^{p}$ of $\Delta_{X}F+ 
 \Delta_{Y}F=0$, that is to say that this function is harmonic in $\R^{2p}$. This is clear if we consider that 
 $F(X+Y)$ is a solution of $\Delta_{X}F= \Delta_{Y}F$.
 
 Remark here that in this analytic continuation, we consider some 
 functions $f(\left|X+\I Y\right|^{2})$, where 
 $$\left|X+\I Y\right|^{2}= \left|X\right|^{2}-\left|Y\right|^{2}+2\I 
 X\cdot Y.$$
 This is not the norm of $X+\I{} Y$ considered as a 
 point in $\C^{p}$.
 
Then,
$F(X+\I\pi(Y))$ is harmonic in $\R^{p+q}$, since the projection of the 
Laplace operator on $\R^{q}$ is the Laplace operator on $\R^{p}$. 
Hence, for any element $R\in SO(q)$, $F(X+\I\pi(RY))$ is harmonic in 
$\R^{p+q}$ since the Laplace operator on $\R^{q}$ is invariant under 
rotations.

From this, we see that 
$$\tilde{H}(X,Y)= \int_{SO(q)}H(X+ \I\pi (RY),0)\ \nu(dR)$$ is harmonic. Observe also 
that it is bi-radial. It is obviously radial in $Y$, since we 
averaged using the Haar measure on $SO(q)$. To see that it is radial 
in $X$, we just observe that, if $R_{1}\in SO(p)$, one has
$$H(R_{1}X+ \I\pi(Z),0)= H\left(X+\I R_{1}^{-1}\pi(Z),0\right)$$ since $H(X,0)$ is 
radial. Moreover, for any $R_{1}\in SO(p)$, there exists 
$R_{2}\in SO(q)$ such that $R_{1}\pi (Z)= \pi(R_{2}Z)$.

Now, let us remark that we will get \eqref{reprOq2} from \eqref{reprebirad} just by expliciting the latter formula. For that, we write $X= \left|X\right|e_{1}$, where $e_{1}\in 
\sphere^{p-1}$ and 
$\pi(RY)= \left|Y\right|\sqrt{\frac{1+u}{2}}Y_{1}$, where $Y_{1}\in \sphere^{p-1}$.

We know that if $R$ is chosen according to the Haar measure on 
$SO(q)$, the law of $R\frac{Y}{\left|Y\right|}$ is uniform on $\sphere^{q}$, and therefore, 
writing $\pi\big(R\frac{Y}{\left|Y\right|}\big)= \left|Y\right| \sqrt{\frac{1+u}{2}}Y_{1}$, the law of $(x,Y_{1})$ 
is 
$\mu_{q-p,p}(dx)\otimes \sigma_{p-1}(dY_{1})$, as we saw at the end 
of paragraph \ref{JacobiPolyn}. Therefore, if we set $s= e_{1}\cdot 
Y_{1}$, the law of $(u,s)$ is $\mu_{p,q-p}(du)\otimes 
\mu_{p-1,p-1}(ds)$.

Then,
$$\left|X+\I\pi(RY)\right|^{2}= 
\left|X\right|^{2}-\frac{1+u}{2}\left|Y\right|^{2} 
+\I\sqrt{2}\left|X\right|\left|Y\right|\sqrt{1+u}e_{1}\cdot Y_{1}\,,$$
 and 
 \begin{multline*}
{ \int_{SO(q)}g\left((\left|X+\I\pi(RY)\right|^{2}\right) \ \nu(dR) = } \\
 \int_{[-1,1]^{2}}g\bigg(\left|X\right|^{2}-\frac{1+u}{2}\left|Y\right|^{2}
 + \I\sqrt{2}\left|X\right|\left|Y\right|\sqrt{1+u}s\bigg)\ \mu_{p,q-p}(dx)
\mu_{p-1,p-1}(ds).
\end{multline*}

To finish the proof of the first formula \eqref{reprebirad}, we observe that the two members of 
\eqref{reprebirad} coincide on $Y=0$. On the other hand, the explicit 
formulation given in \eqref{reprOq2} shows that if $g$ is analytic, 
then the right-hand side in \eqref{reprebirad} is also analytic in 
$(\left|X\right|^{2}, \left|Y\right|^{2})$. Indeed, if we observe 
that the measure $\mu_{p-1,p-1}(ds)$ is symmetric, all odd powers of
$\left|Y\right|\left|X\right|$ in 
the polynomial extension of 
$$\left(\left|X\right|^{2}-\frac{1+u}{2}\left|Y\right|^{2}
+\I\sqrt{2}\left|X\right|\left|Y\right|\sqrt{1+u}s\right)^{n}$$ will 
disappear through integration. We are therefore left with  a series in 
$(\left|X\right|^{2},\left|Y\right|^{2})$.

It remains to see that two analytic harmonic bi-radial functions 
which coincide on $Y=0$ coincide everywhere. Let 
$h(\left|X\right|^{2}, \left|Y\right|^{2})$ an analytic bi-radial 
harmonic function on $\R^{p+q}$.
The function 
$h$ is a solution of 
$$\left(x\partial^{2}_{x}+ \frac{p}{2}\partial_{x}+y\partial^{2}_{y}+ 
\frac{q}{2}\partial_{y}\right)h=0.$$
Then, we see that if we write the expansion
$$h(x,y)= \sum_{n,m}a_{n,m}x^{n}y^{m},$$ one has 
$$a_{n,m+1}=- a_{n+1,m}\frac{(n+1)(n+p/2)}{(m+1)(m+q/2)}.$$
This shows that as soon as one knows $(a_{n,0})$, one knows $h$. This 
completes the proof of Lemma \ref{repreOq}.\qed
\end{proof}

\begin{proof} (Of Koornwinder's formula \eqref{lem2}, for $p$ and $q$ integers. )

In the case where $p$ and $q$ are non-negative integer, it 
turns out that it once again relies on properties of the harmonic 
functions in the Euclidean space.

First, we lift both members on $\R^{p+q}$ and then we multiply  them  
by $R^{k}$, where $R=\left|X\right|^{2}+ \left|Y\right|^{2}$.

As we have seen before, the function
$$(\left|X\right|^{2}+\left|Y\right|^{2})^{k}
P_{k}\left(\frac{\left|X\right|^{2}-\left|Y\right|^{2}}{\left|X\right|^{2}+\left|Y\right|^{2}}\right)$$
is a bi-radial harmonic function, which is a polynomial in 
$\left|X\right|^{2}$ and $\left|Y\right|^{2}$.

It remains to apply Lemma \ref{repreOq} to conclude the proof.\qed
\end{proof}

If we want to extend the proof of Koornwinder's formula when $p$ and 
$q$ are no longer integers, then we just have to observe that, 
setting $S= \left|X\right|^{2}$ and $T= \left|Y\right|^{2}$, we used 
the fact that the fact that
$$H(S,T)=(S+T)^{k} P_{k}\left(\frac{S-T}{S+T}\right)$$is a solution on 
$(0,\infty)^{2}$ of
$$\left(S\partial_{S}^{2}+ \frac{p}{2}\partial_{S}+ T\partial_{T}^{2}+ 
\frac{q}{2}\partial_{T}\right)H=0,$$
and that  for any analytic function $f$,
the function
$$K(S,T)= \int f\left(S-\frac{1+u}{2}T+
 \I\sqrt{2}\sqrt{ST}\sqrt{1+u}s\right)\ \mu_{p,q-p}(dx)
\mu_{p-1,p-1}(ds)$$ is also a solution of the same equation (but this 
time, one has to compute that by brute force!).

\begin{rks}
\item
The proof of Koornwinder's formula gives a representation, for 
analytic functions, of 
solutions $H(S,T)$ of $LH=0$, where
$$L=S\partial_{S}^{2}+ \frac{p}{2}\partial_{S}+ T\partial_{T}^{2}+ 
\frac{q}{2}\partial_{T},$$ in terms of the boundary values 
$H(S,0)$. This is some kind of Poisson formula. In such a 
formula, one has (at least for bounded functions)
$$H(x,y)= \Esp_{x,y}(H(X_{T},Y_{T})),$$
where $(X_{s},Y_{s})$ is the diffusion with generator $L$, and $T$ 
the hitting time of the boundary.

Here, at least when $q\geq2$, the boundary is polar and the set 
$\{y=0\}$ is never attained. But the representation is given here 
through a complex variable (and of course our functions are 
unbounded).  So it happens ``as if'' the process is willing to 
hit the boundary, provided one allows complex values (we do not know 
which is the meaning of that, of course). But it is certainly 
worth looking for more general integral representations of this type, 
with complex values on polar sets.
\bigskip 

\item
Of course, when $q$ converges to $p$, the measure $\mu_{p,q-p}(ds)$ converges to 
the Dirac mass at the point $1$, and Koornwinder's formula of Lemma 
\ref{lem2} gives 
 for the ultraspherical polynomials
\begin{equation}
\frac{P^{p,p}_{k}(x)}{P^{p,p}_{k}(1)} = \int_{[-1,1]^{2}}\left[ 
x+\I\frac{\sqrt{1-x^{2}}}{2}v\right]^{k} \ 
\mu_{p-1,p-1}(dv).
\end{equation}

One may check directly this formula when $p\in \N$, through a much simpler argument, 
using harmonic functions in $\R^{p+1}$ instead of harmonic functions 
in $\R^{2p}$. This time, one has to extend the polynomial 
$P_{k}^{p,p}$ into 
$$F_{k}(X)= \left|X\right|^{k}P^{p,p}_{k}\left(\frac{X}{\left|X\right|}\cdot e_{1}\right),$$ where $e_{1}$ is 
any point of the unit sphere.
\end{rks}
  
\bibliographystyle{plain}  
\bibliography{biblio}

  \end{document}